\documentclass[reqno]{amsart}
\usepackage{amstext,mathrsfs,graphics,color,hyperref,verbatim}

\providecommand{\doi}[1]{}
\def\eqref#1{(\ref{#1})}
\newcommand{\dd}{\,\mrm{d}}
\newcommand{\mbuu}{\sM_{\mrm{bu}}}
\newcommand{\mrm}[1]{\text{\rm #1}}
\newcommand{\cR}{\mathcal{R}} 
\newcommand{\sB}{\mathscr{B}}
\newcommand{\sK}{\mathscr{K}}
\newcommand{\sH}{\mathscr{H}}
\newcommand{\sM}{\mathscr{M}}
\newcommand{\sD}{\mathscr{D}}
\newcommand{\sS}{\mathscr{S}}
\newcommand{\sF}{\mathscr{F}}
\newcommand{\sFb}{\bar{\sF}}
\newcommand{\sC}{\mathscr{C}}
\newcommand{\sMin}{\sM^{m}}
\newcommand{\sMinp}{\sM'{}^{m}}
\newcommand{\set}[2]{\{#1\mid\,#2\}}
\newcommand{\R}{\mathbb{R}}
\newcommand{\N}{\mathbb{N}}
\newcommand{\Z}{\mathbb{Z}}
\newcommand{\Zb}{\overline{\Z}}
\newcommand{\rmax}{\R_{\max}}
\newcommand{\rmaxb}{\overline{\R}_{\max}}
\newcommand{\rmaxbs}{\overline{\R}_{\max}^{_{\scriptstyle S}}}
\newcommand{\rmaxbss}{\overline{\R}_{\max}^{_{\scriptstyle S\times S}}}
\newcommand{\new}[1]{{\em #1}\index{#1}}
\newcommand{\ind}[1]{\chi_{#1}}
\newcommand{\tr}{\mrm{tr}\,}
\newcommand{\pl}{\flat}

\DeclareMathAlphabet{\mathbbold}{U}{bbold}{m}{n}
\newcommand{\zero}{\mathbbold{0}}
\newcommand{\unit}{\mathbbold{1}}
\newcommand{\abar}{\bar A}
\newtheorem{prop}{Proposition}[section]
\newtheorem{corollary}[prop]{Corollary}
\newtheorem{lemma}[prop]{Lemma}
\newtheorem{assumption}[prop]{Assumption}
\newtheorem{theorem}[prop]{Theorem}
\theoremstyle{definition}
\newtheorem{definition}[prop]{Definition}
\theoremstyle{remark}
\newtheorem{example}[prop]{Example}
\newtheorem{remark}[prop]{Remark}
\begin{document}
\title[The max-plus Martin boundary]{The max-plus Martin boundary}
\author{Marianne Akian}
\address{INRIA, Domaine de Voluceau, 
78153 Le Chesnay C\'edex, France}
\email{marianne.akian@inria.fr}
\author{St\'ephane Gaubert}
\address{INRIA, Domaine de Voluceau, 
78153 Le Chesnay C\'edex, France}
\email{stephane.gaubert@inria.fr}
\author{Cormac Walsh}
\address{INRIA, Domaine de Voluceau, 
78153 Le Chesnay C\'edex, France}
\email{cormac.walsh@inria.fr}
\date{June 7, 2005.
}
\thanks{This work was started during a post-doctoral stay of the third author
at INRIA, supported by an ERCIM-INRIA fellowship}

\begin{abstract}
We develop an idempotent version of probabilistic potential theory.
The goal is to describe the set of max-plus harmonic functions,
which give the stationary solutions of deterministic optimal control problems
with additive reward.
The analogue of the Martin compactification is seen to be
a generalisation of the compactification of metric
spaces using (generalised) Busemann functions.
We define an analogue of the minimal Martin boundary and show that it can
be identified with the set of limits of ``almost-geodesics'',
and also the set of (normalised) harmonic functions
that are extremal in the max-plus sense. 
Our main result is a max-plus analogue of the Martin representation
theorem, which represents harmonic functions
by measures supported on the minimal Martin boundary.
We illustrate it by computing the eigenvectors of a class
of translation invariant Lax-Oleinik semigroups.
In this case, we relate the extremal eigenvectors to 
the Busemann points of a normed space.
\end{abstract}
\maketitle
\section{Introduction}
 There exists a correspondence between classical and idempotent
analysis, which was
 brought to light by Maslov and his
 collaborators~\cite{maslov73,maslov92,maslovkololtsov95,litvinov00}.
This correspondence transforms the heat equation
 to an Hamilton-Jacobi equation, and Markov
operators to dynamic programming operators.
So, it is natural
to consider the analogues in idempotent analysis
of harmonic functions, which are
the solutions of the following equation
\begin{align}
\label{sS}
u_i = \sup_{j\in S} (A_{ij} + u_j)
\qquad\mbox{for all $i\in S$.}
\end{align}
The set
 $S$ and the map $A:S\times S\to \R\cup\{-\infty\},\; (i,j)\mapsto A_{ij}$,
which plays the role of the Markov kernel,
are given, and one looks for solutions $u:S\to\R\cup\{-\infty\},\; i\mapsto u_i$.
This equation is the dynamic programming
equation of a deterministic optimal control problem with infinite horizon.
In this context, $S$ is the set of states, the map $A$ gives
the weights or rewards obtained on passing from one state
to another, 
and one is interested  in finding infinite paths
that maximise the sum of the rewards.
Equation~\eqref{sS} is linear
in the max-plus algebra, which is the set $\R\cup\{-\infty\}$
equipped with the operations of maximum
and addition. The term idempotent analysis
refers to the study of structures such as this,
in which the first operation is idempotent.

In potential theory, one uses
the Martin boundary to describe the set of harmonic and
super-harmonic functions of a Markov process,
and the final behaviour of its paths. Our goal here is to obtain
analogous results for Equation~\eqref{sS}.

The original setting for the Martin boundary was
classical potential theory \cite{martin},
where it was used to describe the set of positive solutions of Laplace's
equation.
Doob~\cite{doob} gave a probabilistic interpretation in terms of Wiener
processes and also an extension to the case when time is discrete.
His method was to first establish an integral representation for
super-harmonic functions and then to derive information about
final behaviour of paths. Hunt~\cite{hunt} showed that one could
also take the opposite approach: establish the results concerning
paths probabilistically and then deduce the integral representation.
The approach taken in the present paper is closest to that
of Dynkin~\cite{dynkin}, which contains a simplified version of Hunt's method.

There is a third approach to this subject, using Choquet theory.
However, at present, the tools in the max-plus setting,
are not yet sufficiently developed to allow us to take this route.

Our starting point is the max-plus analogue of the \new{Green kernel},
\begin{equation*}
A^*_{ij}:= \sup\set{A_{i_0i_1}+\dots+A_{i_{n-1}i_{n}}}{n\in\N,\;
  i_0,\dots,i_{n}\in S,\; i_0=i,\; i_n=j} \enspace.
\end{equation*}
Thus, $A^*_{ij}$ is the maximal weight of a path from $i$ to $j$.
We fix a map $i\mapsto \sigma_i$, from $S$ to $\R\cup\{-\infty\}$, which
will play the role of the \new{reference measure}.
We set $\pi_j:=\sup_{k\in S}\sigma_k+A^*_{kj}$.
We define the \new{max-plus Martin space} $\sM$ to be the closure of
the set of maps $\sK:=\set{A^*_{\cdot j}-\pi_{j}}{j\in S}$
in the product topology, and the \new{Martin boundary} to be
$\sM\setminus\sK$. This term must be used with caution
however, since $\sK$ may not be open in $\sM$ (see Example~\ref{ex-triangle}).
The reference measure is often chosen to be a
max-plus Dirac function, taking the value
$0$ at some \new{basepoint} $b\in S$
and the value $-\infty$ elsewhere. In this case,
$\pi_j=A^*_{bj}$.

One may consider the analogue of an ``almost sure'' event to be a set
of outcomes (in our case paths) for which the maximum reward over the
complement is $-\infty$. So we are lead to the notion of
an ``almost-geodesic'', a path of finite total reward, 
see Section~\ref{sec-geo}.
The almost sure convergence of paths in the probabilistic case
then translates into the convergence of every almost-geodesic
to a point on the boundary.

The spectral measure of probabilistic potential theory
also has a natural analogue,
and we use it to give a representation of
the analogues of harmonic functions,
the solutions of~\eqref{sS}.
Just as in probabilistic potential theory,
one does not need the entire Martin boundary
for this representation, a particular subset, called the 
\new{minimal Martin space}, will do.
The probabilistic version is defined in~\cite{dynkin}
to be the set of boundary points for which
the spectral measure is a Dirac measure located at the point itself.
Our definition (see Section~\ref{sec-minmartin})
is closer to an equivalent definition given in the same paper
in which the spectral measure is required only to have a unit of mass at
the point in question. The two definitions are not
equivalent in the max-plus setting and this
is related to the main difference between
the two theories: the representing max-plus measure may not be unique.

Our main theorem (Theorem~\ref{poisson-martin})
is that every (max-plus) harmonic vector $u$ that is integrable
with respect to $\pi$, meaning that $\sup_{j\in S} \pi_j +u_j<\infty$,
can be represented as
\begin{align}\label{e-martin}
u=\sup_{w\in \sMin}\nu(w)+ w,
\end{align}
where $\nu$
is an upper semicontinuous map from the minimal Martin space $\sMin$
to $\R\cup\{-\infty\}$, bounded above.
The map $\nu$ is the analogue of the density of the spectral measure.

We also show that the (max-plus) minimal Martin space is exactly
the set of (normalised) harmonic functions that are
\new{extremal} in the max-plus sense, see Theorem~\ref{th-mr-ext2}. 
We show that each element of the minimal Martin space is either recurrent,
or a boundary point which is the limit of an almost-geodesic
(see Corollary~\ref{cor-geo} and Proposition~\ref{lemma-mr-geo}).

To give a first application of our results,
we obtain in Corollary~\ref{cor-a-tight}
an existence theorem for non-zero harmonic functions of max-plus linear kernels
satisfying a tightness condition, from which we derive a characterisation
of the spectrum of some of these kernels (Corollary~\ref{cor-a-irred}). 

To give a second application, we obtain in Section~\ref{sec-lax}
a representation of the eigenvectors of the Lax-Oleinik semigroup~\cite[\S 3.3]{evans}:
\[
T^tu(x)=\sup_{y\in \R^n} -tL\Big(\frac{y-x}t\Big)+u(y) \enspace,
\]
where $L$ is a convex Lagrangian. This is the evolution semigroup
of the Hamilton-Jacobi equation
\[
\frac{\partial u}{\partial t} = L^\star(\nabla u) \enspace ,
\]
where $L^\star$ denotes the Legendre-Fenchel transform of $L$.
An eigenvector with eigenvalue $\lambda\in\R$
is a function $u$ such that $T^tu=\lambda t+u$
holds for all $t>0$. We compute the eigenvectors
for a subclass of possibly nondifferentiable
Lagrangians (Corollary~\ref{cor-norm1} and Theorem~\ref{cor-final}).

Results and ideas related to the ones of present paper have appeared in several
works: we now discuss them.

Max-plus harmonic functions
have been much studied in the finite-dimensional setting.
The representation formula~\eqref{e-martin}
extends the representation of harmonic vectors
given in the case when $S$ is finite in terms
of the \new{critical} and \new{saturation} graphs.
This was obtained by several authors,
including Romanovski~\cite{romanovski},
Gondran and Minoux~\cite{gondran77} and
Cuninghame-Green~\cite[Th.~24.9]{cuning79}.
The reader may
consult~\cite{maslov92,bcoq,bapat98,gondran02,conv,AGW-s} for more background
on max-plus spectral theory. Relations between max-plus spectral
theory and infinite horizon optimisation are discussed by Yakovenko
and Kontorer~\cite{yakovenko} and Kolokoltsov and Maslov~\cite[\S~2.4]{maslovkololtsov95}. The idea of ``almost-geodesic'' appears there in relation
with ``Turnpike'' theorems.

The max-plus Martin boundary generalises to some extent
the boundary of a metric space defined in terms of (generalised)
Busemann functions by Gromov in~\cite{gromov78}
in the following way
(see also~\cite{gromov} and~\cite[Ch.~II]{ballmann}).
(Note that this is not the same as the Gromov boundary of hyperbolic
spaces.) 
If $(S,d)$ is a complete metric space,
one considers, for all $y,x\in S$, the function $b_{y,x}$ given by
\[
b_{y,x}(z) =
d(x,z)-d(x,y)
\quad
\mrm{for}\; z\in S \enspace .
\]
One can fix the {\em basepoint} $y$ in an arbitrary way.
The space $\sC(S)$ can be equipped
with the topology of uniform convergence on bounded sets,
as in~\cite{gromov78,ballmann}, or with 
the topology of uniform convergence
on compact sets, as in~\cite{gromov}.
The limits of sequences of functions $b_{y,x_n}\in \sC(S)$,
where $x_n$ is a sequence of elements of $S$ going to infinity, 
are called (generalised) \new{Busemann functions}.

When the metric space $S$ is proper,
meaning that all closed bounded subsets of $S$ are compact,
the set of Busemann functions coincides
with the max-plus Martin boundary
obtained by taking
$A_{zx}=A^*_{zx}=-d(z,x)$, and $\sigma$ the
max-plus Dirac function at the basepoint $y$.
This follows from Ascoli's theorem,
see Remark~\ref{rk-loccompact} for details.
Note that our setting is more general since $-A^*$ need not have
the properties of a metric, apart from the triangle inequality
(the case when $A^*$ is not symmetrical is needed
in optimal control).

We note that
Ballman has drawn attention in~\cite[Ch.~II]{ballmann}
to the analogy between this boundary and the probabilistic
Martin boundary. 

The same boundary has recently appeared in the work of Rieffel~\cite{rieffel},
who called it the \new{metric boundary}. 
Rieffel used the term {\em Busemann point} to designate
those points of the metric boundary
that are limits of what he calls ``almost-geodesics''.
We shall see in Corollary~\ref{cor-rieffel}
that these are exactly the points of the
max-plus minimal Martin boundary, at least when $S$ is a proper
metric space. We also relate
Busemann points to extremal eigenvectors of
Lax-Oleinik semigroups, in Section~\ref{sec-lax}.
Rieffel asked in what cases are all boundary points
Busemann points. This problem, as well as the relation between
the metric boundary and other boundaries, has been
studied by Webster and Winchester~\cite{webster,webster03b}
and by Andreev~\cite{andreev}. However, representation problems
like the one dealt with in Theorem~\ref{poisson-martin} do not seem
to have been treated in the metric space context.

Results similar to those of max-plus spectral theory
have recently appeared in weak-KAM theory.
In this context, $S$ is a Riemannian manifold
and the kernel $A$ is replaced by a Lax-Oleinik semigroup,
that is, the evolution semigroup of a Hamilton-Jacobi equation.
Max-plus harmonic functions correspond
to the \new{weak-KAM solutions} of Fathi~\cite{fathi97b,fathi97c,fathi03},
which are essentially the eigenvectors of the Lax-Oleinik semigroup, 
or equivalently,
the viscosity solutions of the ergodic Hamilton-Jacobi equation,
see~\cite[Chapter~8]{fathi03}.
In weak-KAM
theory, the analogue of the Green kernel is called the {\em Ma\~{n}e potential},
the role of the critical graph
is played by the \new{Mather set}, 
and the \new{Aubry set} is related to the saturation graph.
In the case when the manifold is compact,
Contreras~\cite[Theorem~0.2]{contreras}
and Fathi~\cite[Theorem~8.6.1]{fathi03}
gave a representation of the weak-KAM solutions,
involving a supremum of fundamental solutions associated to
elements of the Aubry set.
The case of non-compact manifolds was considered
by Contreras, who defined an analogue of the minimal max-plus Martin
boundary in terms of Busemann functions,
and obtained in~\cite[Theorem~0.5]{contreras}
a representation formula for weak-KAM solutions analogous
to~\eqref{e-martin}. 
Busemann functions also appear in~\cite{fathi03b}.
Other results of weak-KAM theory
concerning non-compact manifolds have been obtained
by Fathi and Maderna~\cite{fathi02}. See also
Fathi and Siconolfi~\cite{fathi04}.
Extremality properties of the elements of
the max-plus Martin boundary (Theorems~\ref{th-mr-ext}
and~\ref{th-mr-ext2} below) do not seem to have been considered
in weak-KAM theory.

Despite the general analogy, the proofs
of our representation theorem for harmonic
functions (Theorem~\ref{poisson-martin}) and
of the corresponding theorems in \cite{contreras} and \cite{fathi03}
require different techniques.
In order
to relate both settings, it would be natural to set $A=B_s$,
where $(B_t)_{t\geq 0}$ is the Lax-Oleinik semigroup, and $s>0$ is
arbitrary.
However, only special kernels $A$ can be written in this way,
in particular $A$ must have an ``infinite divisibility'' property.
Also, not every harmonic function of $B_s$
is a weak-KAM solution associated to the semigroup $(B_t)_{t\geq 0}$.
Thus, the discrete time case is in some sense
more general than the continuous-time case,
but eigenvectors are more constrained in continuous time,
so both settings require distinct treatments.
Nevertheless, in some special cases,
a representation of weak-KAM solutions
follows from our results. This happens
for example in Section~\ref{sec-lax}, where our
assumptions imply that the minimal Martin space
of $B_s$ is independent of $s$.
We note that the Lagrangian there
is not necessarily differentiable,
a property which is required in~\cite{fathi03} and~\cite{contreras}.

The lack of uniqueness of the representing measure is examined
in a further work~\cite{walsh}, where it is shown that
the set of (max-plus) measures representing a given
(max-plus) harmonic function has a least element.

We note that the main results of the present paper have
been announced in the final section of a companion paper,~\cite{AGW-s},
in which max-plus spectral theory was developed under some
tightness conditions.
Here, we use tightness only in Section~\ref{sec-tight}.

\medskip\noindent{\em Acknowledgements.}\/
We thank Albert Fathi for helpful comments, and
in particular for having pointed out to us the work 
of Contreras~\cite{contreras}.
We also thank Arnaud de la Fortelle for references
on the probabilistic Martin boundary theory.

\section{The max-plus Martin kernel and max-plus Martin space}
To show the analogy between the boundary theory
of deterministic optimal control problems
and classical potential theory, it will be convenient
to use max-plus notation.
The \new{max-plus semiring}, $\rmax$, is the set
$\R\cup\{-\infty\}$ equipped with the addition
$(a,b)\mapsto a\oplus b:=\max(a,b)$ and the 
multiplication $(a,b)\mapsto a\odot b:=a+b$. 
We denote by $\zero:=-\infty$ and $\unit:=0$ the zero and unit elements,
respectively. We shall often write $ab$ instead of $a\odot b$.
Since the supremum of an infinite set may be infinite,
we shall occasionally need to consider the {\em completed
max-plus semiring} $\rmaxb$, obtained
by adjoining to $\rmax$ an element $+\infty$,
with the convention that $\zero=-\infty$
remains absorbing for the semiring multiplication.
                                                                                
The sums and products of matrices and vectors are defined in the natural way.
These operators will be denoted by $\oplus$ and concatenation, respectively.
For instance, if $A\in\rmaxbss$, $(i,j)\mapsto A_{ij}$, denotes a matrix
(or kernel), and if $u\in\rmaxbs$, $i\mapsto u_i$ denotes a vector,
we denote by $Au\in\rmaxbs$, $i\mapsto (Au)_i$, the vector defined by
\begin{align*}
(Au)_i:=\bigoplus_{j\in S} A_{ij} u_j \enspace ,
\end{align*}
where the symbol $\oplus$ denotes the usual supremum.

We now introduce the max-plus analogue of the \new{potential kernel}
(Green kernel). Given any matrix $A\in \rmaxbss$, we
define
\begin{align*}
A^*&=I\oplus A\oplus A^2\oplus\cdots \in \rmaxbss \enspace,
\\
A^+&=A\oplus A^2\oplus A^3 \oplus \cdots\in \rmaxbss
\end{align*}
where $I=A^0$ denotes the max-plus identity matrix,
and $A^k$ denotes the $k$th power of the matrix $A$.
The following formulae are obvious:
\begin{align*}
A^*= I \oplus A^+,\qquad A^+=AA^*=A^*A,
\qquad
\mrm{and} \qquad A^*=A^*A^*\enspace .
\end{align*}
It may be useful to keep in mind the graph
representation of matrices:
to any matrix $A\in \rmaxbss$ is associated 
a directed graph with set of nodes $S$ and an arc from $i$ to $j$
if the weight $A_{ij}$ is different from $\zero$.
The weight of a path is by definition the max-plus product (that is, the sum)
of the weights of its arcs. Then, $A^+_{ij}$ and $A^*_{ij}$
represent the supremum of the weights
of all paths from $i$ to $j$ that are, respectively,
of positive an nonnegative length.

Motivated by the analogy with potential theory,
we will say that a vector $u\in\rmax^S$ is (max-plus) \new{harmonic}
if $Au=u$ and \new{super-harmonic} if $Au\leq u$. Note that
we require the entries of a harmonic or super-harmonic
vector to be distinct from $+\infty$.
We shall say that a vector $\pi\in \rmax^S$ is left (max-plus) harmonic
if $\pi A =\pi$,
$\pi$ being thought of as a row vector.
Likewise, we shall say that $\pi$ is left (max-plus) super-harmonic
if $\pi A\leq \pi $.
Super-harmonic vectors have the following elementary characterisation.
\begin{prop}\label{superharm}
A vector $u\in \rmax^S$ is super-harmonic if and only if $u=A^*u$.
\end{prop}
\begin{proof}
If $u\in \rmax^S$ is super-harmonic, then 
$A^ku\leq u$ for all $k\geq 1$, from which it follows that
$u=A^*u$.
The converse also holds, since $AA^*u=A^+u\leq A^*u$.
\end{proof}

{From} now on, we make the following assumption.
\begin{assumption}\label{assump}
There exists a left super-harmonic vector with full support,
in other words a row vector $\pi\in \R^S$ such that $\pi\geq \pi A$.
\end{assumption}
By applying Proposition~\ref{superharm} to the transpose of $A$,
we conclude that $\pi=\pi A^*$.
Since $\pi$ has no components equal to $\zero$,
we see that one consequence of the
above assumption is that  $A^*_{ij}\in\rmax$ for all $i,j\in S$.
A fortiori, $A_{ij}\in \rmax$ for all $i,j\in S$.

The choice of $\pi$ we make will determine which set of harmonic vectors
is the focus of attention. It will be the set of
harmonic vectors $u$ that are \new{$\pi$-integrable},
meaning that $\pi u<\infty$.
Of course, the boundary that we define will also depend on $\pi$, in general.
For brevity, we shall omit the explicit dependence on $\pi$
of the quantities that we introduce and shall omit the assumption
on $\pi$ in the statements of the theorems.
We denote by $\sH$ and $\sS$, respectively, the set
of $\pi$-integrable harmonic and $\pi$-integrable super-harmonic vectors.

It is often convenient to choose
$\pi:=A^*_{b\cdot}$ for some $b\in S$.
(We use the notation $M_{i\cdot}$ and $M_{\cdot i}$ to denote,
respectively, the $i$th row and $i$th column of any matrix $M$.)
We shall say that $b$ is a \new{basepoint} when the vector
$\pi$ defined in this way has finite entries
(in particular, a basepoint has access to every node in $S$).
With this choice of $\pi$, every super-harmonic vector $u\in\rmax^S$ 
is automatically $\pi$-integrable since, by Proposition~\ref{superharm},
$\pi u=(A^* u)_b=u_b<+\infty$.
So, in this case, $\sH$ coincides with the set of 
all harmonic vectors.
This conclusion remains true when $\pi:=\sigma A^*$,
where $\sigma$ is any row vector with finite support,
that is, with $\sigma_i=\zero$ except for finitely many $i$.

We define the \new{Martin kernel} $K$ with respect to $\pi$:
\begin{equation}\label{defi-K}
K_{ij} := A^*_{ij}(\pi_{j})^{-1} \quad \mbox{for all $i,j \in S$} \enspace .
\end{equation}
Since $\pi_{i}A^*_{ij}\leq (\pi A^*)_{j}=\pi_{j}$, we have 
\begin{align}\label{ine-K}
K_{ij}\leq (\pi_{i})^{-1} \quad \mbox{for all $i,j \in S$} \enspace.
\end{align}
This shows that the columns $K_{\cdot j}$ are bounded above
independently of $j$.
By Tychonoff's theorem, the set of columns
$\sK:= \set{K_{\cdot j}}{j\in S}$ is relatively compact 
in the product topology of $\rmax^S$. The \new{Martin space}
$\sM$ is defined to be the closure of $\sK$. 
We call $\sB:=\sM\setminus \sK$ the \new{Martin boundary}.
{From}~\eqref{defi-K} and \eqref{ine-K}, we get that $Aw\leq w$
and $\pi w\leq \unit$ for all $w\in \sK$.
Since the set of vectors with these two properties can be written
\[
\set{w\in\rmax^S}
{\mbox{$A_{ij}w_j\leq w_i$ and $\pi_k w_k\leq \unit$ for all $i,j,k\in S$}}
\]
and this set is obviously closed in the product topology of
$\rmax^S$, we have that
\begin{equation}\label{martin-super}
\sM\subset\sS\quad \text{and} \quad
\pi w\leq \unit \quad \mbox{for all $w\in \sM$}  \enspace .
\end{equation}
\section{Harmonic vectors arising from recurrent nodes}
Of particular interest are those column vectors of $K$ that are harmonic.
To investigate these we will need some basic notions
and facts from max-plus spectral theory.
Define the \new{maximal circuit mean} of $A$ to be
\begin{align*}
 \rho(A) & := \bigoplus_{k\geq 1} (\tr A^k)^{1/k}
\enspace ,
\end{align*}
where $\tr A= \bigoplus_{i\in S} A_{ii}$.
Thus, $\rho(A)$ is the maximum weight-to-length ratio
for all the circuits of the graph of $A$.
The existence of a super-harmonic row vector
with full support, Assumption~\ref{assump}, implies that $\rho(A)\leq \unit$
(see for instance 
Prop.~3.5 of ~\cite{dudnikov} or Lemma~2.2 of~\cite{AGW-s}).
Define the \new{normalised matrix} $\tilde{A}=\rho(A)^{-1} A$.
The max-plus analogue of the notion of recurrence
is defined in~\cite{AGW-s}:
\begin{definition}[Recurrence]
We shall say that a node $i$ is \new{recurrent}
if $\tilde{A}^+_{ii}=\unit$. We denote by
$N^r(A)$ the set of recurrent nodes.
We call \new{recurrent classes} of $A$ the
equivalence classes of $N^r(A)$ with the relation 
$\cR$ defined by $i\cR j$
if $\tilde{A}^+_{ij}\tilde{A}^+_{ji}=\unit$.
\end{definition}
This should be compared with the
definition of recurrence for Markov chains,
where a node is recurrent if one
returns to it with probability one.
Here, a node is recurrent
if we can return to it with reward
$\unit$ in $\tilde{A}$. 

Since $AA^*=A^+\leq A^*$, every column of $A^*$
is super-harmonic. Only those columns
of $A^*$ corresponding to recurrent nodes yield
harmonic vectors:
\begin{prop}[See {\cite[Prop.~5.1]{AGW-s}}]\label{l2}
The column vector $A^*_{\cdot i} $
is harmonic if and only if $\rho(A)=\unit$ and $i$ is recurrent.\qed
\end{prop}
The same is true for the columns of $K$ since they are proportional
in the max-plus sense to those if $A^*$.

The following two results show that 
it makes sense to identify elements in the same
recurrence class.
\begin{prop}\label{egal-K}
Let $i,j\in S$ be distinct.
Then $K_{\cdot i}=K_{\cdot j}$ if and only if
$\rho(A)=\unit$ and $i$ and $j$ are in the same recurrence class.
\end{prop}
\begin{proof}
Let $i,j\in S$ be such that $K_{\cdot i}=K_{\cdot j}$.
Then, in particular, $K_{ii}=K_{ij}$, and so
$A^*_{ij}= \pi_j (\pi_i)^{-1}$. Symmetrically, we obtain 
$A^*_{ji}= \pi_i (\pi_j)^{-1}$. Therefore, $A^*_{ij} A^*_{ji}=\unit$.
If $i\neq j$, then this implies that $A^+_{ii}\geq A^+_{ij} A^+_{ji}=
A^*_{ij} A^*_{ji}=\unit$,
in which case $\rho(A)=\unit$,
$i$ is recurrent, and $i$ and $j$ are in the same
recurrence class. This shows the ``only if'' part of the proposition.
Now let $\rho(A)=\unit$ and $i$ and $j$ be in the same 
recurrence class. Then, according to~\cite[Prop.~5.2]{AGW-s},
$A^*_{\cdot i}= A^*_{\cdot j} A^*_{ji}$,
and so $K_{\cdot i}= K_{\cdot j} (\pi_i)^{-1} \pi_j A^*_{ji}$.
But since $\pi=\pi A^*$, we have that $\pi_i\geq \pi_j A^*_{ji}$,
and therefore $K_{\cdot i}\leq  K_{\cdot j}$. 
The reverse inequality follows from a symmetrical argument.
\end{proof}

\begin{prop}\label{baruconst}
Assume that $\rho(A)=\unit$. Then, for all $u\in \sS$
and $i,\, j$ in the same recurrence class, we have 
$\pi_i u_i=\pi_j u_j$.
\end{prop}
\begin{proof}
Since $\pi\in\R^S$, we can consider the vector
$\pi^{-1}:= (\pi_i^{-1})_{i\in S}$.
That $\pi$ is super-harmonic can be expressed as $\pi_j \geq \pi_i A_{ij}$,
for all $i,j\in S$.
This is equivalent to $(\pi_i)^{-1} \geq A_{ij} (\pi_j)^{-1}$;
in other words, that $\pi^{-1}$, seen as a column vector, is super-harmonic.
Proposition 5.5 of~\cite{AGW-s} states that
the restriction of any two $\rho(A)$-super-eigenvectors
of $A$ to any recurrence class of $A$ are proportional.
Therefore, either $u=\zero$ or
the restrictions of $u$ and $\pi^{-1}$ to any recurrence
class are proportional. In either case, the map $i\in S\mapsto\pi_i u_i$
is constant on each recurrence class.
\end{proof}

\begin{remark}\label{baruconst2}
It follows from these two propositions
that, for any $u\in\sS$, the map $S\to\rmax,\; i\mapsto \pi_i u_i$
induces a map $\sK\to\rmax,\; K_{\cdot i} \mapsto \pi_i u_i$.
Thus, a super-harmonic vector may be regarded as a function defined on $\sK$.
\end{remark}

Let $u\in \rmax^S$ be a $\pi$-integrable vector.
We define the map $\mu_u:\sM\to \rmax$ by
\begin{align*}
\mu_u(w) := \limsup_{K_{\cdot j} \to w} \pi_j {u}_j
 := \inf_{W\ni w} \sup_{K_{\cdot j} \in W} \pi_j {u}_j 
\quad\mrm{for $w\in \sM$} \enspace,
\end{align*}
where the infimum is taken over all neighbourhoods $W$ of $w$ in $\sM$. 
The reason why the limsup above cannot take the value $+\infty$ is
that $\pi_j u_j \leq \pi u <+\infty$ for all $j\in S$.
The following result shows that $\mu_u:\sM\to \rmax$
is an upper semicontinuous extension of the map from $\sK$ to $\rmax$
introduced in Remark~\ref{baruconst2}.
\begin{lemma}\label{u-usc}
Let $u$ be a $\pi$-integrable super-harmonic vector.
Then,
$\mu_u(K_{\cdot i}) = \pi_i u_i$ for each $i\in S$
and
$\mu_u(w)w\le u$ for each $w\in \sM$.
Moreover,
\begin{align*}
 u=\bigoplus_{w\in \sK}  \mu_u(w)  w  
=\bigoplus_{w\in \sM}  \mu_u(w)  w \enspace .
\end{align*}
\end{lemma}
\begin{proof}
By Proposition~\ref{superharm}, $A^* u=u$.
Hence, for all $i\in S$,
\begin{align}\label{uieq}
 u_i &=\bigoplus_{j\in S}  A^*_{ij} u_j= \bigoplus_{j\in S} K_{ij} \pi_j u_j
\enspace .
\end{align}
We conclude that
$u_i\geq K_{ij} \pi_j u_j $ for all $i,j\in S$.
By taking the limsup with respect to $j$ of this inequality,
we obtain that
\begin{equation}\label{pois22}
 u_i\geq \limsup_{K_{\cdot j}\to w} K_{ij} \pi_j u_j 
\geq \liminf_{K_{\cdot j}\to w} K_{ij}  \limsup_{K_{\cdot j}\to w} \pi_j u_j 
= w_i \mu_u(w)\enspace ,\end{equation}
for all $w\in \sM$ and $i\in S$.
This shows the second part of the first assertion of the lemma.
To prove the first part, we apply this inequality with $w=K_{\cdot i}$.
We get that $u_i\geq K_{ii} \mu_u(K_{\cdot i})$.
Since $K_{ii}=(\pi_i)^{-1}$,
we see that $\pi_i u_i \ge \mu_u(K_{\cdot i})$.
The reverse inequality follows from the definition of $\mu_u$.
The final statement of the lemma follows from Equation~\eqref{uieq} and
the first statement.
\end{proof}

\section{The minimal Martin space}\label{sec-minmartin}
In probabilistic potential theory, one does not need the
entire boundary to be able to represent harmonic vectors,
a certain subset suffices.
We shall see that the situation in the max-plus setting is similar.
To define the (max-plus) minimal Martin space, we need to introduce another
kernel:
\begin{align*}
K^{\pl}_{ij} := A^+_{ij}(\pi_{j})^{-1} \quad \mbox{for all $i,j\in S$}\enspace .
\end{align*}
Note that $K^{\pl}_{\cdot j}=A K_{\cdot j}$ is a function of $K_{\cdot j}$.
For all $w\in \sM$, 
we also define $w^{\pl}\in \rmax^S$:
\begin{align*}
w^{\pl}_i=\liminf_{K_{\cdot j}\to w} K^{\pl}_{ij}
\quad \mbox{for all $i\in S$}\enspace .
\end{align*}
The following lemma shows that no ambiguity arises from this
notation since $(K_{\cdot j})^{\pl}=K^{\pl}_{\cdot j}$.
\begin{lemma}\label{w+-prop} 
We have $w^{\pl}=w$ for  $w\in \sB$, and
$w^{\pl}=K^{\pl}_{\cdot j}=A w$ for $w=K_{\cdot j}\in \sK$.
For all $w\in \sM$, we have $w^{\pl}\in\sS$ and $\pi w^{\pl}\leq \unit$.
\end{lemma}
\begin{proof}
Let $w\in \sB$.
Then, for each $i\in S$, there exists a neighbourhood $W$ of $w$ such that
$K_{\cdot i}\not\in W$.
So
\begin{equation*}
w^{\pl}_i
   =\liminf_{K_{\cdot j}\to w} K^{\pl}_{ij}
   =\liminf_{K_{\cdot j}\to w} K_{ij}
   =w_i
\enspace ,
\end{equation*}
proving that $w^{\pl}=w$.

Now let $w=K_{\cdot j}$ for some $j\in S$.
Taking the sequence with constant value
$K_{\cdot j}$, we see that $w^{\pl}\leq K^{\pl}_{\cdot j}$.
To establish the opposite inequality, we observe that
\begin{equation*}
w^{\pl}
   =\liminf_{K_{\cdot k}\to w} A K_{\cdot k}
   \geq \liminf_{K_{\cdot k}\to w} A_{\cdot i} K_{ik}
   =A_{\cdot i} w_i
\qquad
\mbox{for all $i\in S$}
\enspace ,
\end{equation*}
or, in other words, $w^{\pl}\geq Aw$.
Therefore we have shown 
that $w^{\pl}=K^{\pl}_{\cdot j}$.

The last assertion of the lemma
follows from~\eqref{martin-super} and the fact that $\pi$ is super-harmonic.
\end{proof}

Next, we define two kernels $H$ and $H^{\pl}$
over $\sM$.
\begin{align*}
H(z,w):= &\mu_w( z)=\limsup_{K_{\cdot i}\to z} \pi_i{w}_i =
\limsup_{K_{\cdot i}\to z} \lim_{K_{\cdot j}\to w} \pi_i K_{ij}
\\
H^{\pl}(z,w):= &\mu_{w^{\pl}}(z)=\limsup_{K_{\cdot i}\to z} 
\pi_i w^{\pl}_i =
\limsup_{K_{\cdot i}\to z} \, \liminf_{K_{\cdot j}\to w} \pi_i K^{\pl}_{ij}
\enspace .
\end{align*}
Using the fact that $K^{\pl}\leq K$ and Inequality~\eqref{ine-K},
we get that
\begin{align*}
H^{\pl}(z,w)\leq H(z,w)\leq \unit \quad \mbox{for all $w,z\in \sM$}
\enspace.
\end{align*}
If $w\in \sM$, then both $w$ and $w^{\pl}$
are elements of $\sS$ by~\eqref{martin-super} and Lemma~\ref{w+-prop}.
Using the first assertion in Lemma~\ref{u-usc}, we get that
\begin{align}
H(K_{\cdot i},w)&=\pi_i w_i\label{prop0}
\\
H^{\pl}(K_{\cdot i},w)&=\pi_i w^{\pl}_i \enspace .
\end{align}
In particular
\begin{align}
H(K_{\cdot i},K_{\cdot j})&=\pi_i K_{ij}=\pi_i A^*_{ij} (\pi_j)^{-1}
\label{H-prop3}\\
H^{\pl}(K_{\cdot i},K_{\cdot j})&=\pi_i K^{\pl}_{ij}=\pi_i A^+_{ij} (\pi_j)^{-1}
\enspace .\label{H-prop4}
\end{align}
Therefore, up to a diagonal similarity, $H$ and $H^{\pl}$ are
extensions to $\sM\times\sM$ of the kernels $A^*$ and $A^+$
respectively.
\begin{lemma}\label{HH+}
For all $w, z\in \sM$, we have 
\begin{align*}
H(z,w) = \begin{cases}
H^{\pl}(z,w) & \mrm{when } w\neq z\mrm{ or } w=z\in \sB \enspace,\\
\unit & \mrm{otherwise} \enspace . \end{cases}
\end{align*}
\end{lemma}
\begin{proof}
If $w\in \sB$, then $w^\pl=w$ by Lemma~\ref{w+-prop}, and the equality
of $H(z,w)$ and $H^{\pl}(z,w)$ for all $z\in \sM$ follows immediately.

Let $w=K_{\cdot j}$ for some $j\in S$ and let $z\in \sM$
be different from $w$.
Then, there exists a neighbourhood $W$ of $z$
that does not contain $w$.
Applying Lemma~\ref{w+-prop} again, we get that
$w^\pl_i=K^\pl_{ij}=K_{ij}=w_i$ for all $i\in W$.
We deduce that 
$H(z,w)=H^{\pl}(z,w)$ in this case also.

In the final case,
we have $w=z\in\sK$. The result follows from Equation~\eqref{H-prop3}.
\end{proof}

We define the \new{minimal Martin space} to be
\begin{align*}
 \sMin:= \set{w\in \sM}{H^{\pl}(w,w)=\unit}\enspace .
\end{align*}
{From} Lemma~\ref{HH+}, we see that
\begin{align}\label{prop-mr}
\set{w\in \sM}{H(w,w)=\unit }= \sMin\cup \sK \enspace .
\end{align}
\begin{lemma}\label{mrnonzero}
Every $w\in \sMin\cup \sK$ satisfies $\pi w=\unit$.
\end{lemma}
\begin{proof}
We have
\begin{equation*}
\pi {w}
   =   \sup_{i\in S} \pi_i{w}_i
   \ge \limsup_{K_{\cdot i}\to{w}} \pi_i{w}_i
   = H({w},{w})
   = \unit.
\end{equation*}
By Equation~\eqref{martin-super},  $\pi{w}\leq \unit$,
and the result follows.
\end{proof}
\begin{prop}\label{mrsubh} 
Every element of $\sMin$ is harmonic.
\end{prop} 
\begin{proof} 
If $\sK\cap\sMin$ contains an element $w$,
then, from Equation~\eqref{H-prop4}, we see that $\rho(A)=\unit$
and $w$ is recurrent.
It follows from Proposition~\ref{l2} that $w$ is harmonic.

It remains to prove that the same is true for
each element $w$ of $\sB\cap\sMin$.
Let $i\in S$ be such that $w_i\neq \zero$
and assume that $\beta>\unit$ is given.
Since $w\in\sB$, $w$ and $K_{\cdot i}$ will be different.
We make two more observations. Firstly,
by Lemma~\ref{HH+},
$\limsup_{K_{\cdot j}\to w}\pi_j w_j=\unit$.
Secondly, $\lim_{K_{\cdot j}\to w}K_{ij}=w_i$.
{From} these facts, we conclude that there exists $j\in S$,
different from $i$, such that
\begin{equation}
\label{harm1}
\unit \leq \beta \pi_j w_j
\qquad\mbox{and}\qquad
w_i\le \beta K_{ij}
\enspace.
\end{equation}

Now, since $i$ and $j$ are distinct,
we have $A^*_{ij}=A^+_{ij}=(AA^*)_{ij}$.
Therefore, we can find $k\in S$ such that
\begin{equation}
\label{harm2}
A^*_{ij}\le \beta A_{ik} A^*_{kj}
\enspace.
\end{equation}
The final ingredient is that $A^*_{kj} w_j\le w_k$
because $w$ is super-harmonic.
{From} this and the inequalities in~\eqref{harm1} and ~\eqref{harm2},
we deduce that $w_i\le \beta^3 A_{ik} w_k \le \beta^3 (Aw)_i$.
Both $\beta$ and $i$ are arbitrary, so $w\le Aw$.
The reverse inequality is also true since every element of $\sM$ is
super-harmonic. Therefore $w$ is harmonic.
\end{proof}
\section{Martin spaces constructed from different basepoints}\label{sec-mutually}
We shall see that when the left super-harmonic
vector $\pi$ is of the special form $\pi=A^*_{b\cdot}$ for some
basepoint $b\in S$, the corresponding Martin boundary
is independent of the basepoint.
\begin{prop}\label{prop-mutual}
The Martin spaces corresponding to different
basepoints are homeomorphic. 
The same is true for Martin boundaries
and minimal Martin spaces. 
\end{prop}
\begin{proof}
Let $\sM$ and $\sM'$
denote the Martin spaces
corresponding respectively to two different basepoints,
$b$ and $b'$.
We set $\pi=A^*_{b\cdot}$ and
$\pi'=A^*_{b'\cdot}$. 
We denote by $K$ and $K'$ the Martin kernels
corresponding respectively to $\pi$ and $\pi'$.
By construction, $K_{bj}=\unit$ 
holds for all $j\in S$.
It follows that $w_b=\unit$ for all $w\in \sM$.
Using the inclusion in~\eqref{martin-super}, we
conclude that $\sM\subset \sS_b:=\set{w\in \sS}{w_b=\unit}$,
where $\sS$ denotes the set of $\pi$-integrable
super-harmonic functions. 
Observe that $A^*_{bi}$ and $A^*_{b'j}$ are finite for
all $i,j\in S$, since both $b$ and $b'$ are basepoints.
Due to the inequalities $\pi'\geq A^*_{b'b}\pi$
and $\pi\geq A^*_{bb'}\pi'$, 
$\pi$-integrability is equivalent to $\pi'$-integrability.
We deduce that $\sM'\subset \sS_{b'}:=\set{w'\in \sS}{w'_{b'}=\unit}$.
Consider now the maps $\phi$ and $\psi$ defined by
\[
\phi(w)=w(w_{b'})^{-1},\;\forall w\in \sS_b \qquad
 \psi(w')=w'(w'_{b})^{-1},\;\forall w'\in \sS_{b'} \enspace .
\]
Observe that if $w\in \sS_b$, then $w_{b'}\geq A^*_{b'b}w_b=A^*_{b'b}\neq\zero$.
Hence, $w\mapsto w_{b'}$ does not take the value
$\zero$ on $\sS_b$. By symmetry, $w'\mapsto w'_b$ does not take the value zero
on $\sS_{b'}$. It follows that $\phi$
and $\psi$ are mutually inverse homeomorphisms which exchange $\sS_b$
and $\sS_{b'}$. Since $\phi$ sends $K_{\cdot j}$ to $K'_{\cdot j}$,
$\phi$ sends the the Martin space $\sM$,
which is the closure of $\sK:=\set{K_{\cdot j}}{j\in S}$,
to the Martin space $\sM'$, which is the closure
of $\sK':=\set{K'_{\cdot j}}{j\in S}$.
Hence, $\phi$ sends
the Martin boundary $\sM\setminus \sK$ to the Martin boundary
$\sM'\setminus \sK'$.

It remains to show that the minimal Martin space corresponding to $\pi$,
$\sMin$, is sent by $\phi$ to the minimal Martin space corresponding to $\pi'$,
$\sMinp$. Let
\begin{align*}
H'{}^\pl(z',w')&=\limsup_{K'_{\cdot i}\to z'} \, \liminf_{K'_{\cdot j}\to w'} 
A^*_{b'i}A_{ij}^+(A^*_{b'j})^{-1} \enspace .
\end{align*}
Since $\phi$ is an homeomorphism sending $K_{\cdot i}$ to $K'_{\cdot i}$,
a net $(K_{\cdot i})_{i\in I}$ converges to $w$
if and only if the net $(K'_{\cdot i})_{i\in I}$
converges to $\phi(w)$, 
and so
\[
H'{}^\pl(\phi(z),\phi(w))=\limsup_{K_{\cdot i}\to z} \, \liminf_{K_{\cdot j}\to w} A^*_{b'i}A^+_{ij}(A^*_{b'j})^{-1}= 
z_{b'}w_{b'}^{-1} H^\pl(z,w) \enspace .
\]
It follows that $H^\pl(w,w)=\unit$ if and only
if $H'{}^{\pl}(\phi(w),\phi(w))=\unit$.
Hence, $\phi(\sMin)=\sMinp$.
\end{proof}
\begin{remark}
Consider the kernel obtained by symmetrising the kernel $H^\pl$,
\[
(z,w) \mapsto  H^\pl(z,w)H^\pl(w,z) \enspace.
\]
The final argument in the proof of Proposition~\ref{prop-mutual}
shows that this symmetrised kernel is independent
of the basepoint, up to the identification of $w$ and $\phi(w)$.
The same is true for the 
kernel obtained by symmetrising $H$,
\[
(z,w) \mapsto H(z,w)H(w,z)  \enspace .\qedhere
\]
\end{remark}
\section{Martin representation of super-harmonic vectors}\label{sec-super}
In probabilistic potential theory, each super-harmonic vector has a unique
representation as integral over a certain set of vectors, the analogue of
$\sMin\cup \sK$.
The situation is somewhat different in the max-plus setting.
Firstly, according to Lemma~\ref{u-usc}, one does not need the
whole of $\sMin\cup \sK$ to obtain a representation: any set containing
$\sK$ will do.
Secondly, the representation will not necessarily be unique.
The following two theorems, however, show that $\sMin\cup \sK$
still plays an important role.
\begin{theorem}[Martin representation of super-harmonic vectors]\label{Sequal}
For each $u\in \sS$,
$\mu_u$ is the maximal $\nu:\sMin\cup\sK\to\rmax$ satisfying
\begin{align}
u=\bigoplus_{w\in \sMin\cup \sK} \nu(w) w \enspace ,
\label{Sequal1}
\end{align}
Any $\nu:\sMin\cup \sK\to\rmax$ satisfying this equation also satisfies
\begin{align}
\sup_{w\in \sMin\cup\sK} \nu(w) <+\infty \enspace
\label{Sequal1b}
\end{align}
and any $\nu$ satisfying~\eqref{Sequal1b}
defines by~\eqref{Sequal1} an element $u$ of $\sS$.
\end{theorem}
\begin{proof}
By Lemma~\ref{u-usc}, $u$ can be written as~\eqref{Sequal1} with $\nu=\mu_u$.
Suppose that $\nu:\sMin\cup\sK\to\rmax$ is an arbitrary function
satisfying~\eqref{Sequal1}. We have
\begin{equation*}
\pi u=\bigoplus_{w\in \sMin\cup\sK} \nu(w) \pi w \enspace.
\end{equation*}
By Lemma~\ref{mrnonzero},
$\pi w=\unit$ for each $w\in\sMin\cup\sK$.
Since $\pi u<+\infty$, we deduce that~\eqref{Sequal1b} holds.

Suppose that $\nu:\sMin\cup\sK\to\rmax$ is an arbitrary function
satisfying~\eqref{Sequal1b} and define $u$ by~\eqref{Sequal1}.
Since the operation of multiplication by $A$ commutes with arbitrary suprema,
we have $Au\le u$.
Also $\pi u = \bigoplus_{w\in \sMin\cup\sK} \nu(w)<+\infty$.
So $u\in\sS$.

Let $w\in\sMin\cup\sK$.
Then $\nu(w) w_i\leq u_i$ for all $i\in S$.
So we have
\[
\nu(w) H(w,w)
   =    \nu(w) \limsup_{K_{\cdot i}\to w} \pi_i w_i
   \le  \limsup_{K_{\cdot i}\to w} \pi_i u_i
   =    \mu_u(w)\enspace .
\]
Since $H(w,w)=\unit$, we obtain $\nu(w)\leq \mu_u(w)$.
\end{proof}

We shall now give another interpretation of the set $\sMin\cup\sK$.
Let $V$ be a \new{subsemimodule} of $\rmax^S$, that is
a subset of $\rmax^S$ stable under pointwise maximum and
the addition of a constant
(see~\cite{litvinov00,cgq02} for definitions and properties
of semimodules).
We say that a vector
$\xi\in V\setminus\{\zero\}$ is an \new{extremal generator} of $V$ if
$\xi=u\oplus v$ with $u,v\in V$
implies that either $\xi=u$ or $\xi=v$.
This concept has, of course, an analogue in the usual algebra, where
extremal generators are defined for cones.
Max-plus extremal generators are also called \new{join irreducible} elements
in the lattice literature.
Clearly, if $\xi$ is an extremal generator of $V$ then
so is $\alpha \xi$ for all $\alpha\in\R$.
We say that a vector $u\in \rmax^S$
is \new{normalised} if $\pi u=\unit$.
If $V$ is a subset of the set of $\pi$-integrable
vectors, then the set of its extremal generators is
exactly the set of $\alpha \xi$, where $\alpha\in \R$ and
$\xi$ is a normalised extremal generator.
\begin{theorem}\label{th-mr-ext}
The normalised extremal generators of $\sS$ are precisely the elements
of $\sMin\cup \sK$.
\end{theorem}
The proof of this theorem relies on a series of auxiliary results.
\begin{lemma}\label{prop-xi}\sloppy
Suppose that $\xi\in \sMin\cup \sK$ can be written in the form
$\xi=\bigoplus_{w\in \sM} \nu(w) w$,
where $\nu:\sM\to\rmax$ is upper semicontinuous.
Then, there exists $w\in \sM$ such that $\xi=\nu(w) w$.
\end{lemma}
\begin{proof}
For all $i\in S$, we have $\xi_i=\bigoplus_{w\in \sM} \nu(w) w_i$.
As the conventional sum of two upper semicontinuous functions,
the function
$\sM\to\rmax:w \mapsto \nu(w) w_i$ is upper semicontinuous.
Since $\sM$ is compact, the supremum of $\nu(w) w_i$ is attained at
some $w^{(i)}\in \sM$, in other words
$\xi_i=\nu(w^{(i)}) w^{(i)}_i$.
Since $H(\xi,\xi)=\unit$, by definition of $H$,
there exists a net $(i_k)_{k\in D}$ of elements of $S$ such that
$K_{\cdot i_k}$ converges to $\xi$ and 
$\pi_{i_k}\xi_{i_k}$ converges to~$\unit$.
The Martin space $\sM$ is compact and so,
by taking a subnet if necessary,
we may assume that $(w^{(i_k)})_{k\in D}$
converges to some $w\in \sM$.
Now, for all $j\in S$,
\begin{align*}
K_{j i_k}\pi_{i_k}\xi_{i_k}
   &= A^*_{j i_k}\xi_{i_k} 
   =    A^*_{j i_k}\nu(w^{(i_k)}) w^{(i_k)}_{i_k} 
   \le \nu(w^{(i_k)}) w^{(i_k)}_j\enspace ,
\end{align*}
since $w^{(i_k)}$ is super-harmonic.
Taking the limsup as $k\to\infty$, we get that
$\xi_j\le\nu(w)w_j$.
The reverse inequality is true by assumption and therefore
$\xi_j=\nu(w)w_j$.
\end{proof}
The following consequence of this lemma
proves one part of Theorem~\ref{th-mr-ext}.
\begin{corollary}\label{mr-ext}
Every element of $\sMin\cup \sK$ is a normalised extremal generator
of $\sS$.
\end{corollary}
\begin{proof}
Let $\xi\in \sMin\cup \sK$.
We know from Lemma~\ref{mrnonzero} that $\xi$ is normalised.
In particular, $\xi\neq\zero$.
We also know from Equation~\eqref{martin-super} that $\xi\in \sS$.
Suppose $u,v\in \sS$ are such that 
$\xi=u\oplus v$. By Lemma~\ref{u-usc}, we have 
$ u= \bigoplus_{w\in \sM}  \mu_u(w)  w $
and $ v= \bigoplus_{w\in \sM}  \mu_v(w)  w $.
Therefore, $\xi= \bigoplus_{w\in \sM} \nu(w) w$,
with $\nu=\mu_u \oplus \mu_v$. Since $\mu_u$ and 
$\mu_v$ are upper semicontinuous maps from $\sM$ to $\rmax$, so is $\nu$.
By the previous lemma, there exists $w\in \sM$ such
that $\xi= \nu(w) w$.
Now, $\nu(w)$ must equal either $\mu_u(w)$
or $\mu_v(w)$. Without loss of generality, assume the first case.
Then $\xi = \mu_u(w) w\leq u$, and since $\xi\geq u$, we deduce that
$\xi=u$. This shows that $\xi$ is an extremal generator of $\sS$. 
\end{proof} 
The following lemma will allow
us to complete the proof of Theorem~\ref{th-mr-ext}.
\begin{lemma} \label{lemma-ext}
Let $\sF\subset\rmax^S$ have compact closure $\sFb$
in the product topology.
Denote by $V$ the set whose elements are of the form
\begin{align}
\label{e-def-v}
\xi =\bigoplus_{w\in\sF} \nu(w) w\in\rmax^S,
\qquad\mrm{with }
\nu:\sF\to \rmax,\;
\sup_{w\in \sF} \nu(w)<\infty\enspace .
\end{align}
Let $\xi$ be an extremal generator of $V$, and
$\nu$ be as in~\eqref{e-def-v}.
Then, there exists $w\in\sFb$ such that $\xi=\hat{\nu} (w) w$,
where
\begin{align*}
\hat{\nu}(w) :=\limsup_{w'\to w,\, w'\in \sF} \nu(w').
\end{align*}
\end{lemma}
\begin{proof} 
Since $\nu\leq \hat \nu$, we have
$\xi\leq \bigoplus_{w\in\sF} \hat{\nu}(w) w\leq
 \bigoplus_{w\in\sFb} \hat{\nu}(w) w$.
Clearly, $\nu(w) w_i\leq \xi_i$ for all $i\in S$ and $w\in \sF$.
Taking the limsup as $w\to w'$ for any $w'\in\sFb$, we get that
\begin{align*}
\xi_i \geq \hat{\nu}(w')w'_i.
\end{align*}
Combined with the previous inequality, this gives us the representations
\begin{equation}\label{xiegal}
\xi= \bigoplus_{w\in\sF} \hat{\nu}(w) w= 
 \bigoplus_{w\in\sFb} \hat{\nu}(w) w\enspace .
\end{equation}

Consider now, for each $i\in S$ and $\alpha<\unit$,
the set 
\[U_{i,\alpha}:= \set{w\in \sFb}{
 \hat{\nu}(w) w_i
<  \alpha \xi_i}\enspace,
\]
which is open in $\sFb$ since 
the map $w\mapsto \hat{\nu}(w) w_i$ is upper semicontinuous.
Let $\xi\in V\backslash\{\zero\}$ be such that
$\xi\neq \hat{\nu}(w) w$ for all $w\in \sFb$.
We conclude that there exist $i\in S$ and
$\alpha<\unit$ such that $\alpha  \xi_i> \hat{\nu}(w) w_i$,
which shows that $(U_{i,\alpha})_{i\in S,\alpha<\unit}$ is
an open covering of $\sFb$.
Since $\sFb$ is compact, there exists a finite sub-covering
$U_{i_1,\alpha_1},\ldots, U_{i_n,\alpha_n}$.

Using~\eqref{xiegal} and the idempotency of the $\oplus$ law, we
get 
\begin{align}
\xi= \xi^1\oplus \cdots \oplus \xi^n
\qquad\mrm{with }
\xi^j=\bigoplus_{w\in U_{i_j,\alpha_j}\cap\sF} \hat{\nu}(w) w
\enspace,
\label{e-def-xij}
\end{align}
for  $j=1\ldots, n$. 
Since the supremum of $\hat{\nu}$ over $\sFb$ is the same as that over
$\sF$, the vectors
$\xi^1,\ldots,\xi^n$ all belong to $V$. 
Since $\xi$ is an extremal generator of $\sS$, we must
have $\xi=\xi^j$ for some $j$. Then $U_{i_j,\alpha_j}\cap\sF$ is non-empty,
and so $\xi_{i_j}>\zero$.
But,
from the definition of $U_{i_j,\alpha_{j}}$,
\begin{align*}
\xi^j_{i_j}=\bigoplus_{w\in U_{i_j,\alpha_j}\cap \sF} 
\hat{\nu}(w) w_{i_j}
&\leq \alpha_{i_j}\xi_{i_j}
   < \xi_{i_j}
\enspace.
\end{align*}
This shows that $\xi^j$ is different from $\xi$, and so
Equation~\eqref{e-def-xij} gives the required decomposition of
$\xi$, proving it is not an extremal generator of $V$.
\end{proof}
We now conclude the proof of Theorem~\ref{th-mr-ext}:
\begin{corollary}
\label{cor-extrofs}
Every normalised extremal generator of $\sS$ belongs
to $\sMin\cup \sK$.
\end{corollary}
\begin{proof}
Take $\sF=\sMin\cup \sK$ and let $V$ be as defined in Lemma~\ref{lemma-ext}.
Then, by definition, $\sFb=\sM$, which is compact.
By Theorem~\ref{Sequal}, $V=\sS$. 
Let $\xi$ be a normalised
extremal generator of $\sS$. Again by Theorem~\ref{Sequal},
$\xi=\oplus_{w\in \sF} \mu_{\xi}(w) w$. 
Since $\mu_{\xi}$ is upper semicontinuous on $\sM$,
Lemma~\ref{lemma-ext} yields $\xi=\mu_{\xi} (w) w$ for some $w\in \sM$,
with $\mu_{\xi}(w)\neq\zero$ since $\xi\neq\zero$.
Note that $\mu_{\alpha u}=\alpha\mu_u$ for all $\alpha\in\rmax$
and $u\in \sS$. Applying this to the previous equation and evaluating
at $w$, we deduce that  $\mu_{\xi}(w)=\mu_{\xi}(w) \mu_w(w)$.
Thus, $H(w,w)=\mu_w(w)=\unit$.
In addition, $\xi$ is normalised
and so, by Lemma~\ref{mrnonzero}, 
\begin{equation*}
\unit= \pi\xi=\mu_{\xi} (w) \pi w= \mu_{\xi} (w).
\end{equation*}
Hence
$\xi=w \in  \sMin\cup \sK$. 
\end{proof}
\section{Almost-geodesics}\label{sec-geo}
In order to prove a Martin representation
theorem for harmonic vectors,
we will use a notion appearing in~\cite{yakovenko} 
and~\cite[\S~2.4]{maslovkololtsov95}, which
we will call almost-geodesic. A variation
of this notion appeared in~\cite{rieffel}.
We will compare the two notions later in the section.

Let $u$ be a super-harmonic vector,
that is $u\in\rmax^S$ and  $Au\leq u$.
Let $\alpha\in\rmax$ be such that $\alpha\geq\unit$.
We say that a sequence $(i_k)_{k\geq 0}$ with values in $S$
is an \new{$\alpha$-almost-geodesic} with respect to $u$
if $u_{i_0}\in \R$ and
\begin{equation}\label{eqn:usemidef}
u_{i_0} \leq \alpha A_{i_0 i_1} \cdots  A_{i_{k-1}i_k} u_{i_k}
\quad\mbox{for all $k\geq 0\enspace$}.
\end{equation}
Similarly, $(i_k)_{k\geq 0}$ is an \new{$\alpha$-almost-geodesic}
with respect to a left super-harmonic vector $\sigma$ if
$\sigma_{i_0}\in \R$ and
\begin{align*}
\sigma_{i_k} \leq \alpha \sigma_{i_0} A_{i_0 i_1}
\cdots  A_{i_{k-1}i_k} 
\quad\mbox{for all $k\geq 0\enspace$}.
\end{align*}
We will drop the reference to $\alpha$ when its value is unimportant.
Observe that, if $(i_k)_{k\geq 0}$ is an almost-geodesic with respect to 
some right super-harmonic vector $u$, then both
$u_{i_k}$ and $A_{i_{k-1} i_k}$ are in $\R$ for all $k\geq 0$.
This is not necessarily true if $(i_k)_{k\geq 0}$
is an almost-geodesic with respect to a left super-harmonic vector $\sigma$, 
however, if additionally $\sigma_{i_k}\in\R$ for all $k\geq 0$,
then $A_{i_{k-1} i_k}\in\R$ for all $k\geq 0$.

\begin{lemma}\label{geo-ri-le}
Let $u,\sigma\in\rmax^S$ be, respectively, right and left super-harmonic
vectors and assume that $u$ is $\sigma$-integrable, that is
$\sigma u<+\infty$. If $(i_k)_{k\geq 0}$ is an almost-geodesic
with respect to $u$, and if $\sigma_{i_0}\in\R$, then $(i_k)_{k\geq 0}$ 
is an almost-geodesic with respect to $\sigma$.
\end{lemma}
\begin{proof}
Multiplying Equation~(\ref{eqn:usemidef})
by $\sigma_{i_k}(u_{i_0})^{-1}$, we obtain
\begin{align*} \sigma_{i_k}&\leq \alpha \sigma_{i_k} u_{i_k}(u_{i_0})^{-1}
A_{i_0 i_1} \cdots  A_{i_{k-1}i_k}
\leq \alpha (\sigma u )(\sigma_{i_0} u_{i_0})^{-1} \sigma_{i_0}
A_{i_0 i_1} \cdots  A_{i_{k-1}i_k}\enspace.
\end{align*}
So $(i_k)_{k\geq 0}$ 
is a  $\beta$-almost-geodesic with respect to $\sigma$,
with $\beta:=\alpha (\sigma u ) (\sigma_{i_0} u_{i_0})^{-1}\geq \alpha $.
\end{proof}

\begin{lemma}\label{semigeo-beta}
Let $(i_k)_{k\geq 0}$ be an almost-geodesic with respect to $\pi$ and
let $\beta>\unit$. Then, for $\ell$ large enough,
$(i_k)_{k\geq \ell}$ is a $\beta$-almost-geodesic with respect to $\pi$.
\end{lemma}
\begin{proof}
Consider the matrix $\abar_{ij}:=\pi_i A_{ij}(\pi_j)^{-1}$.
The fact that $(i_k)_{k\geq 0}$ is an $\alpha$-almost-geodesic
with respect to $\pi$ is equivalent to
\begin{align*}
p_k:= 
(\abar_{i_{0}i_1})^{-1}
\cdots 
(\abar_{i_{k-1}i_k})^{-1}
\leq \alpha
\quad\mbox{for all $k\geq 0$}
\enspace.
\end{align*}
Since $(\abar_{i_{\ell-1}i_\ell})^{-1}\geq \unit$ for all
$\ell\geq 1$, the sequence $\{p_k\}_{k\geq 1}$ is nondecreasing.
The upper bound then implies it converges to a finite limit.
The Cauchy criterion states that
\[
\lim_{\ell,k \to \infty, \, \ell< k} 
\abar_{i_{\ell}i_{\ell+1}}
\cdots \abar_{i_{k-1}i_k}= \unit\enspace .
\]
This implies that, given any $\beta>\unit$,
$\abar_{i_{\ell}i_{\ell+1}} \cdots \abar_{i_{k-1}i_k}
   \geq \beta^{-1}$ for $k$ and $\ell$ large enough, with $k>\ell$.
Writing this formula in terms of $A$ rather than $\abar$,
we see that, for $\ell$ large enough, 
$(i_k)_{k\geq \ell}$ is a $\beta$-almost-geodesic with respect to $\pi$.
\end{proof}
\begin{prop}\label{semigeo-conv}
If $(i_k)_{k\geq 0}$ is an almost-geodesic with respect to $\pi$, then
$K_{\cdot i_k}$ converges to some $w\in \sMin$.
\end{prop}
\begin{proof}
Let $\beta>\unit$.
By Lemma~\ref{semigeo-beta},
$(i_k)_{k\geq \ell}$ is a $\beta$-almost-geodesic with respect to $\pi$,
for $\ell$ large enough.
Then, for all $k>\ell$,
\begin{align*}
\pi_{i_k}
   \leq \beta \pi_{i_{\ell}} A^+_{i_{\ell} i_k}
   \leq \beta \pi_{i_{\ell}} A^*_{i_{\ell} i_k} 
\enspace.
\end{align*}
Since $\pi$ is left super-harmonic, we have
$\pi_{i_{\ell}} A^*_{i_{\ell} i_k}   \leq \pi_{i_k}$.
Dividing by $\beta\pi_{i_k}$ the former inequalities, we
deduce that
\begin{equation}\label{limKellk}
\beta^{-1} \leq \pi_{i_{\ell}}  K^{\pl}_{i_{\ell} i_k} \leq 
\pi_{i_{\ell}} K_{i_{\ell} i_k} \leq \unit\enspace.
\end{equation}
Since $\sM$ is compact, it suffices to check that all
convergent subnets of $K_{\cdot i_k}$ have
the same limit $w\in \sMin$. Let $(i_{k_d})_{d\in D}$
and $(i_{\ell_e})_{e\in E}$ denote subnets of $(i_k)_{k\geq 0}$,
such that the nets $(K_{\cdot i_{k_d}})_{d\in D}$ and $(K_{\cdot i_{\ell_e}})_{e\in E}$ converge to some $w\in \sM$ and $w'\in \sM$, respectively.
Applying~\eqref{limKellk} with $\ell=\ell_e$ and $k=k_d$, and
taking the limit with respect to $d$,
we obtain $\beta^{-1}\leq 
\pi_{i_{\ell_e}} w_{i_{\ell_e}}$. Taking now 
the limit with respect to $e$,  we get that $\beta^{-1}\leq 
H(w',w)$. Since this holds for all $\beta>\unit$, we obtain
$\unit\leq H(w',w)$, thus $H(w',w)=\unit$.
{From} Lemma~\ref{u-usc},
we deduce that $w\geq \mu_{w}(w')w'
=H(w',w)w'=w'$. By symmetry, we conclude that $w=w'$,
and so $H(w,w)=\unit$.
By Equation~\eqref{prop-mr}, 
$w\in \sMin\cup\sK$.
Hence, $(K_{\cdot i_{k}})_{k\geq 0}$
converges towards some $w\in \sMin\cup \sK$.

Assume by contradiction that $w\not\in \sMin$.
Then, $w=K_{\cdot j}$ for some $j\in S$,
and $H^{\pl}(w,w)<\unit$ by definition
of $\sMin$.
By~\eqref{H-prop4},
this implies that $\pi_j K^{\pl}_{jj}=A^+_{jj}<\unit$.
If the sequence
$(i_k)_{k\geq 0}$ takes the value $j$ infinitely often,
then, we can deduce from Equation~\eqref{limKellk} that $A^+_{jj}=\unit$,
a contradiction.
Hence, for $k$ large enough, $i_k$ does
not take the value $j$, which implies, by Lemma~\ref{w+-prop},
that $w_{i_{k}}=w^{\pl}_{i_{k}}$. 
Using Equation~\eqref{limKellk}, 
we obtain $H^{\pl}(w,w)\geq \limsup_{k\to\infty} 
\pi_{i_{k}} w^{\pl}_{i_{k}}=\limsup_{k\to\infty} 
\pi_{i_{k}} w_{i_{k}}=\unit$, which contradicts our assumption on $w$.
We have shown that $w\in\sMin$.
\end{proof}
\begin{remark}\label{rk-rieffel}
An inspection of the proof of Proposition~\ref{semigeo-conv}
shows that the same conclusion holds under the weaker hypothesis
that for all $\beta>\unit$, we have $\pi_{i_k}\leq \beta\pi_{i_\ell}A^+_{i_\ell i_k}$ for all $\ell$ large enough and $k> \ell$.
\end{remark}
\begin{lemma}\label{lem-tech-geo}
If  $(i_k)_{k\geq 0}$ is an almost-geodesic with respect to $\pi$,
and if $w$ is the limit of $K_{\cdot i_k}$, then
\[
\lim_{k\to\infty} \pi_{i_k}w_{i_k}=\unit \enspace .
\]
\end{lemma}
\begin{proof}
Let $\beta>\unit$. By Lemma~\ref{semigeo-beta},
$(i_k)_{k\geq \ell}$ is a $\beta$-almost-geodesic with respect to $\pi$
for $\ell$ large enough. Hence, for all $k\geq \ell$, 
$\pi_{i_k} \leq \beta \pi_{i_\ell}A^*_{i_\ell i_k}$,
and so $\unit \leq \beta \pi_{i_\ell}A^*_{i_\ell i_k}\pi_{i_k}^{-1}
=\beta \pi_{i_\ell} K_{i_\ell i_k}$. Since $K_{i_\ell i_k}$ converges
to $w_{i_\ell}$ when $k$ tends to infinity, we deduce that $\unit\leq \beta
\liminf_{\ell\to \infty} \pi_{i_\ell}w_{i_\ell}$, and since this holds
for all $\beta>\unit$, we get 
$\unit \leq \liminf_{\ell\to\infty}\pi_{i_\ell}w_{i_\ell}$. Since $\pi_jw_j\leq \unit$ for all $j$, the lemma is proved.
\end{proof}
\begin{prop}\label{prop-carac-specmes}
Let $u$ be a $\pi$-integrable super-harmonic vector. Then, $\mu_u$
is continuous along almost-geodesics, meaning that
if $(i_k)_{k\geq 0}$ is an almost-geodesic with respect to $\pi$
and if $K_{\cdot i_k}$ tends to $w$, then, 
\[
\mu_u(w) = \lim_{k\to\infty} \mu_u(K_{\cdot i_k})= \lim_{k\to\infty} \pi_{i_k}u_{i_k} \enspace .
\]
\end{prop}
\begin{proof}
Recall that $\pi_{i}u_{i}=\mu_u(K_{\cdot i})$ holds for all $i$,
as shown in Lemma~\ref{u-usc}. It also follows
from this lemma that $u\geq \mu_u(w)w$, and so
$\pi_i u_i \geq \pi_i w_i \mu_u(w)$ for all $i\in S$. 
Hence,
\begin{align*}
\liminf_{k\to\infty} \pi_{i_k}u_{i_k}&\geq \liminf_{k\to\infty}
\pi_{i_k}w_{i_k}\mu_u(w)\\
&=\mu_u(w) \enspace,
\end{align*}
by Lemma~\ref{lem-tech-geo}. 
Moreover, $\limsup_{k\to\infty} \pi_{i_k}u_{i_k} \leq \mu_u(w)$,
by definition of $\mu_u(w)$.
\end{proof}
Combining Lemma~\ref{geo-ri-le} and Proposition~\ref{semigeo-conv},
we deduce the following.
\begin{corollary}\label{cor-geo}
If $(i_k)_{k\geq 0}$ is an almost-geodesic with respect to a
$\pi$-integrable super-harmonic vector, then
$K_{\cdot i_k}$ converges to some element of $\sMin$.
\end{corollary}
For brevity, we shall say sometimes
that an almost-geodesic $(i_k)_{k\geq 0}$ \new{converges}
to a vector $w$ when $K_{\cdot i_k}$ converges to $w$.
We state a partial converse to Proposition~\ref{semigeo-conv}.
\begin{prop} \label{lemma-mr-geo} Assume that $\sM$ is
first-countable.
For all $w\in\sMin$, there exists an almost-geodesic with respect to $\pi$
converging to $w$.
\end{prop}
\begin{proof}
By definition, $H^{\pl}(w,w)=0$.
Writing this formula explicitly in terms of $A_{ij}$ and making the
transformation $\abar_{ij}:=\pi_i A_{ij}(\pi_j)^{-1}$,
we get
\begin{align*}
\limsup_{K_{\cdot i}\to w} \liminf_{K_{\cdot j}\to w}
\abar^+_{ij} = \unit\enspace .\end{align*} 
Fix a sequence $(\alpha_k)_{k\geq 0}$ 
in $\rmax$ such that $\alpha_k>\unit$ and
$\alpha:=\alpha_0 \alpha_1\cdots<+\infty$. 
Fix also a decreasing sequence $(W_k)_{k\geq 0}$ of
open neighbourhoods of $w$.
We construct a sequence $(i_k)_{k\geq 0}$ in $S$ inductively as follows.
Given $i_{k-1}$, we choose $i_{k}$ to have the following three
properties:
\begin{enumerate}\renewcommand{\theenumi}{\alph{enumi}}
\item
\label{firstitem}
$K_{\cdot i_{k}}\in W_{k}$,
\item
\label{seconditem}
$\liminf_{K_{\cdot j}\to w}
\abar^+_{i_k j}  > \alpha_{k}^{-1}$,
\item
\label{thirditem}
$\abar^+_{i_{k-1} i_{k}} > \alpha_{k-1}^{-1}$.
\end{enumerate}
Notice that it is possible to satisfy~(\ref{thirditem}) because
$i_{k-1}$ was chosen to satisfy~(\ref{seconditem}) at the previous step.
We require $i_0$ to satisfy~(\ref{firstitem}) and (\ref{seconditem})
but not (\ref{thirditem}).
Since $\sM$ is first-countable, one can choose
the sequence
$(W_k)_{k\geq 0}$ in such a way that every sequence
$(w_k)_{k\geq 0}$ in $\sM$ 
with $w_k\in W_k$ converges to $w$.
By~(\ref{thirditem}), one can find, for all $k\in\N$, a finite sequence
$(i^\ell_k)_{0\leq \ell\leq N_k}$ such that $i^0_k=i_k$, $i^{N_k}_k=i_{k+1}$,
and
\[ \abar_{i_{k}^0, i_{k}^1} \cdots 
\abar_{i_{k}^{N_{k}-1}, i_{k}^{N_k}}  > \alpha_k^{-1}\quad
\mbox{for all $k\in\N$} \enspace .\]
Since $\abar_{ij}\leq \unit$ for all $i,j\in S$,
we obtain 
\begin{align*}
\abar_{i_{k}^0, i_{k}^1} \cdots 
\abar_{i_{k}^{n-1}, i_{k}^{n}}  > \alpha_k^{-1}\quad
\mbox{for all $k\in\N$, $1\leq n\leq N_k$}\enspace .
\end{align*}
Concatenating the sequences $(i^\ell_k)_{0\leq \ell\leq N_k}$,
we obtain a sequence $(j_m)_{m\geq 0}$ such that 
$\alpha^{-1}\leq \abar_{j_{0} j_{1}} \cdots 
\abar_{j_{m-1} j_{m}}$  for all $m\in\N$,
in other words an $\alpha$-almost-geodesic with respect to $\pi$.
{From} Lemma~\ref{semigeo-conv}, we know that $K_{\cdot j_m}$ converges to
some point in $\sM$.
Since $(i_k)$ is a subsequence of $(j_m)$ and
$K_{\cdot i_k}$ converges to $w$, we deduce that 
$K_{\cdot j_m}$ also converges to $w$.
\end{proof}
\begin{remark}
If $S$ is countable, the product topology on $\sM$ is metrisable.
Then, the assumption of Proposition~\ref{lemma-mr-geo} is satisfied.
\end{remark}
\begin{remark}\label{rk-loccompact}
Assume that $(S,d)$ is a metric space, let $A_{ij}=A^*_{ij}=-d(i,j)$
for $i,j\in S$, and let $\pi=A^*_{b\cdot}$ for any $b\in S$.
We have $K_{\cdot j}=-d(\cdot,j)+d(b,j)$.
Using the triangle inequality for $d$, we see that,
for all $k\in S$, the function $K_{\cdot k}$ is non-expansive,
meaning that $|K_{ik}-K_{jk}|\leq d(i,j)$ for all $i,j\in S$. 
It follows that every map in $\sM$ is non-expansive.
By Ascoli's theorem, the topology of pointwise convergence
on $\sM$ coincides with the topology
of uniform convergence on compact sets.
Hence, if $S$ is a countable union of compact sets,
then $\sM$ is metrisable and
the assumption of Proposition~\ref{lemma-mr-geo}
is satisfied. 
\end{remark}
\begin{example}\label{ex-denum}
The assumption in Proposition~\ref{lemma-mr-geo} cannot
be dispensed with.
To see this, take $S=\omega_1$, the first uncountable ordinal.
For all $i,j\in S$,
define $A_{ij}:=0$ if $i<j$ and $A_{ij}:=-1$ otherwise.
Then, $\rho(A)=\unit$ and $A=A^+$.
Also $A^*_{ij}$ equals $0$ when $i\leq j$
and $-1$ otherwise.
We take $\pi:=A^*_{0 \cdot}$,
where $0$ denotes the smallest ordinal.
With this choice, $\pi_i=\unit$ for all $i\in S$, and $K=A^*$.

Let $\sD$ be the set of maps $S\to\{-1,0\}$ that are non-decreasing
and take the value 0 at 0.
For each $z\in\sD$, define 
$s(z):=\sup\set{i\in S}{z_i=0}\in S\cup\{\omega_1\}$.
Our calculations above lead us to conclude that
\begin{equation*}
\sK = \set{z\in \sD}{\mbox{$s(z)\in S$ and $z_{s(z)}=0$}} \enspace .
\end{equation*}
We note that
$\sD$ is closed in the product topology on $\{-1,0\}^S$ and contains $\sK$.
Furthermore, every $z\in\sD\setminus \sK$ is the limit of the net
($A^*_{\cdot d})_{d\in D}$ indexed by the directed set
$D=\set{d\in S}{d<s_z}$. 
Therefore the Martin space is given by 
$\sM=\sD$.
Every limit ordinal $\gamma$ less than or equal to $\omega_1$
yields one point $z^\gamma$ in the Martin boundary $\sB:=\sM\setminus\sK$:
we have $z^\gamma_i=0$ for $i<\gamma$, and $z^\gamma_i=-1$
otherwise.

Since $A^+_{ii}=A_{ii}=-1$ for all $i\in S$,
there are no recurrent points, and so $\sK\cap \sMin$ is empty.
For any $z\in \sB$, we have $z_d=0$ for all $d<s(z)$.
Taking the limsup, we conclude that
$H(z,z)=\unit$, thus $\sMin=\sB$.
In particular, the identically zero vector $z^{\omega_1}$ is in $\sMin$.

Since a countable union of countable sets is countable,
for any sequence $(i_k)_{k\in\N}$ of elements of $S$, 
the supremum $I=\sup_{k\in \N} i_k$ belongs to $S$, 
and so its successor ordinal,
that we denote by $I+1$, also belongs to $S$.
Since $\lim_{k\to\infty} K_{I+1, i_k}=-1$,
$K_{\cdot i_k}$ cannot converge to $z^{\omega_1}$,
which shows that the point $z^{\omega_1}$ in the minimal Martin space
is not the limit of an almost-geodesic.
\end{example}
We now compare our notion of almost-geodesic with that
of Rieffel~\cite{rieffel} in the metric space
case. We assume that
$(S,d)$ is a metric space and take $A_{ij}=A^*_{ij}=-d(i,j)$ and
$\pi_j=-d(b,j)$, for an some $b\in S$.
The compactification of $S$ discussed in~\cite{rieffel},
called there the \new{metric compactification},
is the closure of $\sK$
in the topology of uniform convergence on compact sets,
which, by Remark~\ref{rk-loccompact},
is the same as its closure
in the product topology. It thus
coincides with the Martin space $\sM$.
We warn the reader that variants of the metric
compactification can be found in the literature,
in particular, the references~\cite{gromov78,ballmann} use the topology
of uniform convergence on bounded sets rather than on compacts.

Observe that
the basepoint $b$ can be chosen in an arbitrary way: indeed,
for all $b'\in S$,
setting $\pi'=A^*_{b'\cdot}$, we get
$\pi'\geq A^*_{b'b}\pi$ and $\pi\geq A^*_{bb'}\pi'$,
which implies that almost-geodesics in our sense are the same
for the basepoints $b$ and $b'$. Therefore, 
when speaking of almost-geodesics in our sense, in a metric space,
we will omit the reference to $\pi$.

Rieffel defines an almost-geodesic as an $S$-valued map
$\gamma$ from an unbounded set $\mathcal{T}$ of real nonnegative
numbers containing $0$, such that for all $\epsilon>0$,
for all $s\in \mathcal T$ large enough, and for all $t\in \mathcal T$ such that $t\geq s$,
\[|d(\gamma(t),\gamma(s))+d(\gamma(s),\gamma(0))-t|<\epsilon
\enspace .\]
By taking $t=s$, one sees that $|d(\gamma(t),\gamma(0))-t|<\epsilon$.
Thus, almost-geodesics in the sense of Rieffel are ``almost''
parametrised by arc-length, unlike those in our sense.
\begin{prop}\label{prop-rieffel}
Any almost-geodesic in the sense of Rieffel has a subsequence
that is an almost-geodesic in our sense. Conversely, 
any almost-geodesic in our sense that is not bounded
has a subsequence that is an almost-geodesic in the sense of Rieffel.
\end{prop}
\begin{proof}
Let $\gamma:\mathcal{T} \to S$
denote an almost-geodesic in the sense of Rieffel.
Then, for all $\beta>\unit$, we have
\begin{align}
A^*_{\gamma(0),\gamma(t)}
\leq \beta A^*_{\gamma(0),\gamma(s)}A^*_{\gamma(s)\gamma(t)} \enspace
\label{e-superadd}
\end{align}
for all $s\in \mathcal{T}$ large enough and for all $t\in \mathcal{T}$
such that $t\geq s$. 
Since the choice of the basepoint $b$
is irrelevant, we may assume that
$b=\gamma(0)$,
so that $\pi_{\gamma(s)}=A^*_{\gamma(0),\gamma(s)}$.
As in the proof of Lemma~\ref{semigeo-beta}
we set $\bar A_{ij}=\pi_i A^*_{ij}\pi_j^{-1}$.
We deduce from~\eqref{e-superadd} that
\[
\beta^{-1}\leq \bar A_{\gamma(s)\gamma(t)}\leq \unit
\enspace .
\] 
Let us choose a sequence $\beta_1,\beta_2,\ldots\geq\unit $
such that the product $\beta_1\beta_2\ldots$ converges
to a finite limit.
We can construct a sequence
$t_0<t_1<\ldots$ of elements of $\mathcal{T}$
such that, setting $i_k=\gamma(t_{i_k})$,
\[\bar A_{i_ki_{k+1}} \geq \beta_k^{-1}
\enspace . \]
Then, the product
$\bar A_{i_0i_{1}}
\bar A_{i_1i_{2}}\cdots$
converges, which implies that the sequence $i_0$, $i_1,\ldots$
is an almost-geodesic in our sense.

Conversely, let $i_0,i_1,\ldots$ be an almost-geodesic in our sense, and
assume that $t_k=d(b,i_k)$ is not bounded. After replacing
$i_k$ by a subsequence, we may assume that $t_0<t_1<\ldots$.
We set $\mathcal{T}=\{t_0,t_1,\ldots\}$ and $\gamma(t_k)=i_k$.
We choose the basepoint $b=i_{0}$, so that $t_{0}=0\in \mathcal{T}$,
as required in the definition of Rieffel.
Lemma~\ref{semigeo-beta} implies that
\[ A^*_{bi_k}\leq \beta A^*_{bi_\ell}A^*_{i_\ell i_k} 
\]
holds for all $\ell$ large enough and for all $k\geq \ell$.
Since $t_{k}^{-1}=A^*_{bi_k}$, 
$\gamma$ is an almost-geodesic in the sense of Rieffel.
\end{proof}
Rieffel called the limits of almost-geodesics
in his sense \new{Busemann points}. 
\begin{corollary}\label{cor-rieffel}
Let $S$ be a proper metric space.
Then the minimal Martin space is the disjoint union
of $\sK$ and of the set of Busemann points of $S$.
\end{corollary}
\begin{proof}
Since $A^+_{ii}=-d(i,i)=0$ for all $i$, the set $\sK$ is included
in the minimal Martin space $\sMin$. We next show that $\sMin\setminus \sK$
is the set of Busemann points.

Let $w\in \sM$ be a Busemann point.
By Proposition~\ref{prop-rieffel} we can find
an almost-geodesic in our sense $i_0,i_1,\ldots$ such that 
 $K_{\cdot i_k}$ converges to $w$ and $d(b,i_k)$ is unbounded.
We know from Proposition~\ref{semigeo-conv}
that $w\in\sMin$. It remains to check that $w\not\in \sK$.
To see this, we show that for all $z\in \sM$,
\begin{align}
\lim_{k\to \infty} H(K_{\cdot i_k}, z) = H(w,z) \enspace .
\label{e-cont}
\end{align}
Indeed, for all $\beta>\unit$, 
letting $k$ tend to infinity in~\eqref{limKellk}
and using~\eqref{prop0}, we get
\[ \beta^{-1} \leq 
\pi_{i_{\ell}} w_{i_\ell}=
H(K_{\cdot i_\ell}, w) \leq \unit \enspace,
\]
for $\ell$ large enough.
Hence,
$\lim_{\ell\to\infty} H(K_{\cdot i_\ell}, w) = \unit$.
By Lemma~\ref{u-usc}, $z\geq H(w,z)w$. We deduce
that
 $H(K_{\cdot i_\ell},z) \geq H(w,z)H(K_{\cdot i_\ell},w)$,
and so $\liminf_{\ell\to\infty} H(K_{\cdot i_\ell},z) \geq H(w,z)$.
By definition of $H$, 
$\limsup_{\ell\to\infty} H(K_{\cdot i_\ell},z) \leq 
\limsup_{K_{\cdot j}\to w} H(K_{\cdot j},z)=H(w,z)$,
which shows~\eqref{e-cont}. Assume now that $w\in \sK$,
that is, $w=K_{\cdot j}$ for some $j\in S$,
and let us apply~\eqref{e-cont} to $z=K_{\cdot b}$. 
We have $H(K_{\cdot i_k},z)=A^*_{bi_k}A^*_{i_kb}=-2\times d(b,i_k)\to -\infty$.
Hence, $H(w,z)=-\infty$. But
$H(w,z)=A^*_{bj}A^*_{jb}=-2\times d(b,j)>-\infty$,
which shows that $w\not\in \sK$.

Conversely, let $w\in \sMin\setminus \sK$. By Proposition~\ref{lemma-mr-geo},
$w$ is the limit of an almost-geodesic in our sense. 
Observe that this almost-geodesic is unbounded. Otherwise,
since $S$ is proper,
$i_k$ would have a converging subsequence,
and by continuity of the map $i \mapsto K_{\cdot i}$, we would have $w\in \sK$,
a contradiction. It follows from Proposition~\ref{prop-rieffel} that $w$
is a Busemann point. 
\end{proof}
\section{Martin representation of harmonic vectors}\label{sec-harm}

\begin{theorem}[Poisson-Martin representation of harmonic vectors]
\label{poisson-martin}
Any element $u\in\sH$ can be written as
\begin{align}
u=\bigoplus_{w\in \sMin} \nu(w) w \enspace ,
\label{Hequal1}
\end{align}
with $\nu:\sMin\to\rmax$, and necessarily,
\begin{align*}
\sup_{w\in \sMin} \nu(w) <+\infty \enspace .
\end{align*}
Conversely, any $\nu:\sMin\to\rmax$
satisfying the latter inequality
defines by~\eqref{Hequal1} an element $u$ of $\sH$.
Moreover, given $u\in \sH$,
$\mu_u$ is the maximal $\nu$ satisfying~\eqref{Hequal1}.
\end{theorem}
\begin{proof}
Let $u\in \sH$. Then $u$ is also in $\sS$ and so, from Lemma~\ref{u-usc},
we obtain that 
\begin{equation}\label{e-int10}
u= \bigoplus_{w\in \sM} \mu_u(w)  w \geq \bigoplus_{w\in \sMin}
\mu_u(w) w \enspace .
\end{equation}
To show the opposite inequality, let us fix some 
$i\in S$ such that $u_i\neq \zero$.
Let us also fix some sequence $(\alpha_k)_{k\geq 0}$ 
in $\rmax$ such that $\alpha_k>\unit$ for all $k\ge 0$ and such that
$\alpha:=\alpha_0\alpha_1\cdots<+\infty$. 
Since $u=Au$, one can construct a sequence
$(i_k)_{k\geq 0}$ in $S$ starting at $i_0:=i$, and such that
\begin{align*}
u_{i_k}&\leq \alpha_k A_{i_{k} i_{k+1}} u_{i_{k+1}}\quad 
\mbox{for all $k\geq 0$}\enspace . 
\end{align*}
Then, 
\begin{align}\label{ine-repr}
u_{i_0}&\leq \alpha  A_{i_{0} i_{1}}\cdots A_{i_{k-1} i_{k}} u_{i_{k}}
\leq \alpha  A^*_{i_{0} i_{k}} u_{i_{k}}\quad
\mbox{for all $k\geq 0$}\enspace ,
\end{align}
and so $(i_k)_{k\geq 0}$ is an $\alpha$-almost-geodesic
with respect to $u$. Since $u$ is $\pi$-integrable, 
we deduce using Corollary~\ref{cor-geo} that $K_{\cdot i_k}$ converges
to some $w\in \sMin$.
{From}~\eqref{ine-repr}, we get 
$u_i\leq \alpha K_{i i_k} \pi_{i_k} u_{i_k}$,
and letting $k$ go to infinity,
we obtain $u_i\leq \alpha w_i \mu_u(w)$.
We thus obtain
\[ u_i\leq \alpha \bigoplus_{w\in \sMin} \mu_u(w) w_i \enspace .\]
Since $\alpha$ can be chosen arbitrarily close to $\unit$, we deduce
the inequality opposite to~\eqref{e-int10}, which shows
that~\eqref{Hequal1} holds with $\nu=\mu_u$.

The other parts of the theorem are proved in a manner similar to
Theorem~\ref{Sequal}.
\end{proof}
\begin{remark}
The maximal representing measure $\mu_u$ at every point
that is the limit of an almost geodesic
can be computed by taking the limit of $\pi_i u_i$ along 
any almost-geodesic converging to this point.
See Proposition~\ref{prop-carac-specmes}.
\end{remark}
In particular, $\sH=\{\zero\}$ if and only if $\sMin$ is empty.
We now prove the analogue of Theorem~\ref{th-mr-ext}
for harmonic vectors.
\begin{theorem}\label{th-mr-ext2}
The normalised extremal generators of $\sH$ are precisely the elements
of $\sMin$.
\end{theorem}
\begin{proof}
We know from Theorem~\ref{th-mr-ext} that each element of $\sMin$
is a normalised extremal generator of $\sS$. Since $\sH\subset\sS$, 
and $\sMin\subset \sH$ (by Proposition~\ref{mrsubh}),
this implies that each element of $\sMin$
is a normalised extremal generator of $\sH$. 

Conversely, by the same arguments
as in the proof of Corollary~\ref{cor-extrofs},
taking $\sF=\sMin$ in Lemma~\ref{lemma-ext} and
using Theorem~\ref{poisson-martin} instead of Lemma~\ref{u-usc},
we get that each normalised extremal generator $\xi$ of $\sH$
belongs to $\sMin\cup \sK$. 
Since, by Proposition~\ref{l2},
no element of $\sK\setminus \sMin$ can be harmonic,
we have that $\xi\in \sMin$.
\end{proof}

\begin{remark}\label{rk-finrec}
Consider the situation when there are only finitely many recurrence classes
and only finitely many non-recurrent nodes.
Then $\sK$ is a finite set, so that
$\sB$ is empty, $\sM=\sK$,
and $\sMin$ coincides with the set of columns
$K_{\cdot j}$ with $j$ recurrent. The
representation theorem (Theorem~\ref{poisson-martin})
shows in this case that each harmonic vector is a finite max-plus linear
combination of the recurrent columns of $A^*$, as is the case in
finite dimension.
\end{remark}

\section{Product Martin spaces}\label{sec-decomp}
In this section,
we study the situation where the set $S$ is the Cartesian product of two sets,
$S_1$ and $S_2$, and $A$ and $\pi$ can be decomposed
as follows:
\begin{align}
A=A_1\otimes I_2 \oplus I_1\otimes A_2 \enspace,
\quad 
\pi=\pi_1\otimes \pi_2\enspace .\label{form}
\end{align}
Here, $\otimes$ denotes the max-plus tensor product of matrices or vectors,
$A_i$ is an $S_i\times S_i$ matrix, $\pi_i$ is a vector
indexed by $S_i$,
and $I_i$ denotes the $S_i\times S_i$ max-plus identity matrix. For instance,
$(A_1\otimes I_2)_{(i_1,i_2),(j_1,j_2)}=(A_1)_{i_1j_1}(I_2)_{i_2j_2}$,
which is equal to $(A_1)_{i_1j_1}$ if $i_2=j_2$, and to $\zero$ otherwise.
We shall always assume that $\pi_i$ is left super-harmonic with respect to $A_i$, for $i=1,2$.
We denote by $\sM_i$ the corresponding Martin space, by $K_i$ the corresponding
Martin kernel, etc.

We introduce the map
\[
\imath: \rmax^{S_1}\times \rmax^{S_2}\to \rmax^S ,\; 
\imath(w_1,w_2)=w_1\otimes w_2 \enspace ,
\]
which is obviously continuous for the product topologies.
The restriction of $\imath$ to the set of $(w_1,w_2)$
such that $\pi_1w_1=\pi_2w_2=\unit$ is injective.
Indeed, if $w_1\otimes w_2=w'_1\otimes w'_2$,
applying the operator $I_1\otimes \pi_2 $ on both sides
of the equality, we get $w_1\otimes \pi_2w_2=w'_1\otimes \pi_2w'_2$,
from which we deduce that $w_1=w'_1$ if $\pi_2w_2=\pi_2w'_2=\unit$.

\begin{prop}\label{prop-tensor}
Assume that $A$ and $\pi$ are of the form~\eqref{form},
and that $\pi_iw_i=\unit$ for all $w_i\in \sM_i$ and $i=1,2$.
Then, the map $\imath$
is a homeomorphism from
$\sM_1\times \sM_2$ to the Martin space $\sM$ of $A$,
and sends $\sK_1\times \sK_2$ to $\sK$.
Moreover, the same map sends
\[
\sMin_1 \times (\sK_2\cup  \sMin_2)
\;\cup \;
(\sK_{1}\cup  \sMin_1)\times 
\sMin_2
\]
to the minimal Martin space $\sMin$
of $A$.
\end{prop}
The proof of Proposition~\ref{prop-tensor} relies on several
lemmas.
\begin{lemma}\label{lem-formula}
If $A$ is given by~\eqref{form}, then, $A^*=A_1^*\otimes
A_2^*$ and 
\[A^+=A_1^+\otimes A_2^*\oplus A_1^*\otimes A_2^+\enspace .\]
\end{lemma}
\begin{proof}
Summing the equalities
$A^k=\bigoplus_{1\leq \ell \leq k} A_1^\ell \otimes A_2^{k-\ell}$,
we obtain $A^*=A_1^*\otimes A_2^*$. Hence,
$A^+=AA^*=(A_1\otimes I_2\oplus I_1\otimes A_2)(A_1^*\otimes A_2^*)
=A_1^+\otimes A_2^*\oplus A_1^*\otimes A_2^+$.
\end{proof}
We define the kernel $H\circ\imath$ 
from $(\sM_1\times \sM_2)^2$ to $\rmax$, by 
$H\circ\imath((z_1,z_2),(w_1,w_2))=H(\imath(z_1,z_2),\imath(w_1,w_2))$.
The kernel $H^{\pl}\circ \imath$ is defined
from $H^{\pl}$
in the same way.
\begin{lemma}\label{lem-hom}
If $A^*=A_1^*\otimes A_2^*$ and $\pi=\pi_1\otimes\pi_2$, 
then $\sK=\imath(\sK_1\times \sK_2)$ and $\imath(\sM_1\times \sM_2)=\sM$.
Moreover, if $\pi_iw_i=\unit$ for all $w_i\in \sM_i$ and $i=1,2$, then
$\imath$ is an homeomorphism from $\sM_1\times \sM_2$ to $\sM$,
and $H\circ\imath=H_1\otimes H_2$.
\end{lemma}
\begin{proof}
Observe that $K=K_1\otimes K_2$. Hence, $\sK=\imath(\sK_1\times \sK_2)$.
Let $\overline{X}$ denote the closure of any set $X$.
Since $\overline{\sK_i}=\sM_i$,
we get $\overline{\sK_1\times \sK_2}
=\sM_1\times \sM_2$, and so 
$\overline{\sK_1\times \sK_2}$ is compact.
Since $\imath$ is continuous, we deduce
that $\imath(\overline{\sK_1\times \sK_2})=
\overline{\imath(\sK_1\times \sK_2)}$.
Hence, $\imath(\sM_1\times \sM_2)
=\overline{\sK}=\sM$. Assume now that $\pi_iw_i=\unit$ for all $w_i\in \sM_i$
and $i=1,2$, so that the restriction
of $\imath$ to $\sM_1\times \sM_2$ is injective.
Since $\sM_1\times \sM_2$ is compact, we deduce
that $\imath$ is an homeomorphism
from $\sM_1\times \sM_2$ to its image, that is, $\sM$.
Finally, let $z=\imath(z_1,z_2)$ and $w=\imath(w_1,w_2)$,
with $z_1,w_1\in\sM_1$ and $z_2,w_2\in\sM_2$.
Since $\imath$ is an homeomorphism
from $\sM_1\times \sM_2$ to $\sM$, we can write
$H(z,w)$ in terms of limsup and limit for the product topology
of $\sM_1\times \sM_2$:
\begin{align}
H(z,w)=\limsup_{\scriptstyle (K_1)_{\cdot i_1}\to z_1\atop\scriptstyle
(K_2)_{\cdot i_2}\to z_2}
\lim_{\scriptstyle (K_1)_{\cdot j_1}\to w_1\atop\scriptstyle
(K_2)_{\cdot j_2}\to w_2}
\pi_{(i_1,i_2)}K_{(i_1,i_2),(j_1,j_2)} \enspace .
\label{e-prod}
\end{align}
Since $A^*=A_1^*\otimes A_2^*$ and $\pi=\pi_1\otimes\pi_2$,
we can write the right hand side term of~\eqref{e-prod}
as the product of two terms that are
both bounded from above:
\[\pi_{(i_1,i_2)}K_{(i_1,i_2),(j_1,j_2)}=
\left((\pi_1)_{i_1}(K_1)_{i_1,j_1}\right)\left((\pi_2)_{i_2}(K_2)_{i_2,j_2}\right)
\enspace .\]
Hence, the limit and limsup in~\eqref{e-prod} become a product
of limits and limsups, respectively, and so
$H(z,w)=H_1(z_1,w_1)H_2(z_2,w_2)$.
\end{proof}
\begin{lemma}\label{lemma-comp}
Assume that $A$ and $\pi$ are of the form~\eqref{form} and that 
$\pi_iw_i=\unit$
for all $w_i\in \sM_i$ and $i=1,2$.
Then
\begin{align}
H^\pl\circ \imath&=H_1^\pl\otimes H_2 \oplus H_1\otimes H_2^\pl \enspace .\label{e-hpl}
\end{align}
\end{lemma}
\begin{proof}
By Lemma~\ref{lem-formula}, 
$A^+=A_1^+\otimes A_2^*\oplus A_1^*\otimes A_2^+$,
and so 
\[
K^{\pl}=K_1^\pl\otimes K_2\oplus K_1\otimes K_2^\pl \enspace .
\]
Let $z=\imath(z_1, z_2)$
and $w=\imath(w_1,w_2)$,
with $z_1,w_1\in\sM_1$, $z_2,w_2\in\sM_2$.
In a way similar to~\eqref{e-prod}, we can
write $H^{\pl}$ as
\begin{align*}
H^\pl(z,w)&= \limsup_{\scriptstyle(K_1)_{\cdot i_1}\to z_1\atop\scriptstyle
(K_2)_{\cdot i_2}\to z_2}
\liminf_{\scriptstyle (K_1)_{\cdot j_1}\to w_1\atop\scriptstyle
(K_2)_{\cdot j_2}\to w_2} \pi_{(i_1,i_2)}K^{\pl}_{(i_1,i_2),(j_1,j_2)}
\enspace .
\end{align*}
The right hand side term is a sum of products:
\begin{align*}
\pi_{(i_1,i_2)}K^{\pl}_{(i_1,i_2),(j_1,j_2)}
&=
(\pi_1)_{i_1}(K_1^{\pl})_{i_1j_1}(\pi_2)_{i_2}(K_2)_{i_2j_2}
\oplus 
(\pi_1)_{i_1}(K_1)_{i_1 j_1} 
(\pi_2)_{i_2}(K_2^{\pl})_{i_2j_2} \enspace .
\end{align*}
We now use the following two general observations.
Let $(\alpha_d)_{d\in D}$, $(\beta_e)_{e\in E}$, $(\gamma_d)_{d\in D}$,
$(\delta_e)_{e\in E}$ be nets of elements of $\rmax$ that are bounded from above.
Then, 
\[\limsup_{d,e} \alpha_d\beta_e\oplus \gamma_d \delta_e
=(\limsup_d \alpha_d)(\limsup_e \beta_e) \oplus (\limsup_d\gamma_d)(\limsup_e\delta_e)
\enspace .
\]If additionally the nets
$(\beta_e)_{e\in E}$ and $(\gamma_d)_{d\in D}$ converge, 
we have
\[\liminf_{d,e} \alpha_d\beta_e\oplus \gamma_d \delta_e
=(\liminf_d \alpha_d)(\lim_e \beta_e) \oplus (\lim_d\gamma_d)(\liminf_e\delta_e)\enspace.\]
Using both identities, we deduce that $H^{\pl}$ is given by~\eqref{e-hpl}.
\end{proof}
\begin{proof}[Proof of Proposition~\ref{prop-tensor}] We know
from Lemma~\ref{lem-formula} that $A^*=A_1^*\otimes A_2^*$,
and so, by Lemma~\ref{lem-hom}, $\imath$ is an homeomorphism
from $\sM_1\times \sM_2$ to $\sM$. 
Since the kernels $H_1,H_1^\pl,H_2$ and $H_2^\pl$ all take values
less than or equal to $\unit$, we conclude from~\eqref{e-hpl}
that, when $z=\imath(z_1,z_2)$, $H^\pl(z,z)=\unit$
if and only if $H_1^\pl(z_1,z_1)=H_2(z_2,z_2)=\unit$
or $H_1(z_1,z_1)=H_2^{\pl}(z_2,z_2)=\unit$.
Using Equation~\eqref{prop-mr} and the definition of the minimal Martin space,
we deduce that 
\[\sMin=\imath\big(\sMin_1 \times (\sK_2\cup  \sMin_2)
\;\cup \;
(\sK_{1}\cup  \sMin_1)\times 
\sMin_2\big) \enspace .\qedhere\]
\end{proof}
\begin{remark}
The assumption that $\pi_iw_i=\unit$ for all $w_i\in\sM_i$
is automatically satisfied when the left super-harmonic
vectors $\pi_i$ originate from basepoints,
that is, when $\pi_i=(A_i)^*_{b_i,\cdot}$ for some basepoint
$b_i$. Indeed, we already observed in the proof of
Proposition~\ref{prop-mutual} that every vector $w_i\in\sM_i$
satisfies $(\pi_i)_{b_i}(w_i)_{b_i}=\unit$.
By~\eqref{martin-super}, $\pi_iw_i\leq\unit$.
We deduce that $\pi_iw_i= \unit$.
\end{remark}
\begin{remark}
Rieffel~\cite[Prop.~4.11]{rieffel}
obtained a 
version of the first part of Lemma~\ref{lem-hom}
for metric spaces. His result
states that if $(S_1,d_1)$ and $(S_2,d_2)$ are locally compact
metric spaces, 
and if their product $S$ is equipped with the sum of the metrics,
$d((i_1,i_2),(j_1,j_2))=d_1(i_1,j_1)+d_2(i_2,j_2)$, then the metric
compactification of $S$ can be identified with the Cartesian
product of the metric compactifications of $S_1$ and $S_2$. 
This result can be re-obtained from Lemma~\ref{lem-hom}
by taking $(A_1)_{i_1,j_i}=-d_1(i_1,j_1)$,
$(A_2)_{i_2,j_2}=-d_2(i_2,j_2)$, $\pi_{i_1}=-d_1(i_1,b_1)$,
and $\pi_{i_2}=-d(i_2,b_2)$, for arbitrary
basepoints $b_1,b_2\in\Z$. We shall illustrate this in Example~\ref{ex-3}.
\end{remark}
\section{Examples and Counter-Examples}
We now illustrate our results
and show various features that the Martin space may have.

\begin{example}\label{ex-1}
Let $S=\N$, $A_{i,i+1}=0$ for all $i\in\N$,
$A_{i,0}=-1$ for all $i\in \N\setminus\{0\}$ and 
$A_{ij}=-\infty$ elsewhere. We choose the basepoint $0$,
so that $\pi=A^*_{0,\cdot}$.
The graph of $A$ is:
\begin{center}
\begin{picture}(0,0)%
\includegraphics{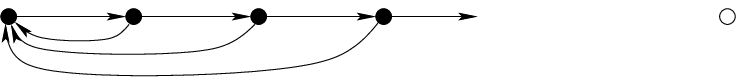}%
\end{picture}%
\setlength{\unitlength}{1973sp}%
\begingroup\makeatletter\ifx\SetFigFont\undefined%
\gdef\SetFigFont#1#2#3#4#5{%
  \reset@font\fontsize{#1}{#2pt}%
  \fontfamily{#3}\fontseries{#4}\fontshape{#5}%
  \selectfont}%
\fi\endgroup%
\begin{picture}(7066,955)(2618,-4240)
\put(4276,-3561){\makebox(0,0)[lb]{\smash{{\SetFigFont{6}{7.2}{\rmdefault}{\mddefault}{\itdefault}{\color[rgb]{0,0,0}$0$}%
}}}}
\put(5551,-3561){\makebox(0,0)[lb]{\smash{{\SetFigFont{6}{7.2}{\rmdefault}{\mddefault}{\itdefault}{\color[rgb]{0,0,0}$0$}%
}}}}
\put(6676,-3561){\makebox(0,0)[lb]{\smash{{\SetFigFont{6}{7.2}{\rmdefault}{\mddefault}{\itdefault}{\color[rgb]{0,0,0}$0$}%
}}}}
\put(5116,-4028){\makebox(0,0)[lb]{\smash{{\SetFigFont{6}{7.2}{\rmdefault}{\mddefault}{\itdefault}{\color[rgb]{0,0,0}$-1$}%
}}}}
\put(4051,-3886){\makebox(0,0)[lb]{\smash{{\SetFigFont{6}{7.2}{\rmdefault}{\mddefault}{\itdefault}{\color[rgb]{0,0,0}$-1$}%
}}}}
\put(3151,-3561){\makebox(0,0)[lb]{\smash{{\SetFigFont{6}{7.2}{\rmdefault}{\mddefault}{\itdefault}{\color[rgb]{0,0,0}$0$}%
}}}}
\put(3083,-3826){\makebox(0,0)[lb]{\smash{{\SetFigFont{6}{7.2}{\rmdefault}{\mddefault}{\itdefault}{\color[rgb]{0,0,0}$-1$}%
}}}}
\end{picture}%

\end{center}
States (elements of $S$) are represented by black dots.
The white circle represents the extremal boundary element $\xi$,
that we next determine.
In this example, $\rho(A)=\unit$, and $A$ has no recurrent class.
We have $A^*_{ij}=\unit$ for $i\leq j$ and $A^*_{ij}=-1$
for $i>j$, so the Martin space of $A$ corresponding
to $\pi=A^*_{0\cdot}$ consists of the columns $A^*_{\cdot j}$,
with $j\in \N$, together with the vector $\xi$ whose entries
are all equal to $\unit$.
We have $\sB=\{\xi\}$. 
One can easily check that $H(\xi,\xi)=\unit$.
Therefore, $\sMin=\{\xi\}$. 
Alternatively, we may use Proposition~\ref{semigeo-conv}
to show that $\xi\in\sMin$,
since $\xi$ is the limit of the almost-geodesic
$0,1,2,\ldots$.
Theorem~\ref{poisson-martin} says that $\xi$ is
the unique (up to a multiplicative constant) non-zero harmonic vector.
\end{example}
\begin{example}\label{ex-2}
Let us modify Example~\ref{ex-1} by setting $A_{00}=0$,
so that the previous graph becomes:
\begin{center}
\begin{picture}(0,0)%
\includegraphics{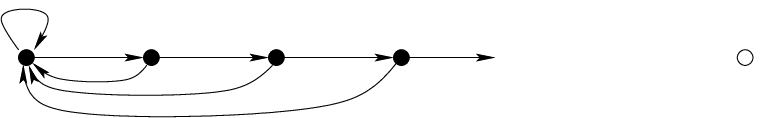}%
\end{picture}%
\setlength{\unitlength}{1973sp}%
\begingroup\makeatletter\ifx\SetFigFont\undefined%
\gdef\SetFigFont#1#2#3#4#5{%
  \reset@font\fontsize{#1}{#2pt}%
  \fontfamily{#3}\fontseries{#4}\fontshape{#5}%
  \selectfont}%
\fi\endgroup%
\begin{picture}(7236,1236)(2448,-4240)
\put(4276,-3561){\makebox(0,0)[lb]{\smash{{\SetFigFont{6}{7.2}{\rmdefault}{\mddefault}{\itdefault}{\color[rgb]{0,0,0}$0$}%
}}}}
\put(5551,-3561){\makebox(0,0)[lb]{\smash{{\SetFigFont{6}{7.2}{\rmdefault}{\mddefault}{\itdefault}{\color[rgb]{0,0,0}$0$}%
}}}}
\put(6676,-3561){\makebox(0,0)[lb]{\smash{{\SetFigFont{6}{7.2}{\rmdefault}{\mddefault}{\itdefault}{\color[rgb]{0,0,0}$0$}%
}}}}
\put(5116,-4028){\makebox(0,0)[lb]{\smash{{\SetFigFont{6}{7.2}{\rmdefault}{\mddefault}{\itdefault}{\color[rgb]{0,0,0}$-1$}%
}}}}
\put(4051,-3886){\makebox(0,0)[lb]{\smash{{\SetFigFont{6}{7.2}{\rmdefault}{\mddefault}{\itdefault}{\color[rgb]{0,0,0}$-1$}%
}}}}
\put(3151,-3561){\makebox(0,0)[lb]{\smash{{\SetFigFont{6}{7.2}{\rmdefault}{\mddefault}{\itdefault}{\color[rgb]{0,0,0}$0$}%
}}}}
\put(3083,-3826){\makebox(0,0)[lb]{\smash{{\SetFigFont{6}{7.2}{\rmdefault}{\mddefault}{\itdefault}{\color[rgb]{0,0,0}$-1$}%
}}}}
\put(2626,-3136){\makebox(0,0)[lb]{\smash{{\SetFigFont{6}{7.2}{\rmdefault}{\mddefault}{\itdefault}{\color[rgb]{0,0,0}$0$}%
}}}}
\end{picture}%

\end{center}
We still have $\rho(A)=\unit$,
the node $0$ becomes recurrent, and the minimal Martin space
is now $\sMin=\{K_{\cdot 0},\xi\}$, where $\xi$ is defined
in Example~\ref{ex-1}.
Theorem~\ref{poisson-martin} says that every harmonic vector
is of the form $\alpha K_{\cdot 0} \oplus \beta\xi $,
that is $\sup(\alpha+ K_{\cdot 0},\beta+\xi)$ with the
notation of classical algebra, for some $\alpha,\beta\in \R\cup\{-\infty\}$.
\end{example}

\begin{example}\label{ex-z}
Let $S=\Z$, $A_{i,i+1}=A_{i+1,i}=-1$ for $i\in \Z$,
and $A_{ij}=\zero$ elsewhere.
We choose $0$ to be the basepoint, so that $\pi=A^*_{0,\cdot}$.
The graph of $A$ is:
\begin{center}
\begin{picture}(0,0)%
\includegraphics{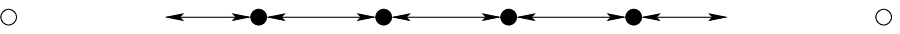}%
\end{picture}%
\setlength{\unitlength}{1973sp}%
\begingroup\makeatletter\ifx\SetFigFont\undefined%
\gdef\SetFigFont#1#2#3#4#5{%
  \reset@font\fontsize{#1}{#2pt}%
  \fontfamily{#3}\fontseries{#4}\fontshape{#5}%
  \selectfont}%
\fi\endgroup%
\begin{picture}(8566,1202)(218,-4187)
\end{picture}%

\end{center}
We are using the same conventions as in the previous examples,
together with the following additional conventions:
the arrows are bidirectional since the matrix is symmetric, 
and each arc has weight $-1$ unless otherwise specified. 
This example and the next were considered by Rieffel~\cite{rieffel}.

We have $\rho(A)=-1<\unit$, which implies there are no recurrent nodes.
We have $A^*_{i,j}=-|i-j|$, and so
$K_{i,j}=|j|-|i-j|$. There are
two Martin boundary points, $\xi^{+}=\lim_{j\to\infty} K_{\cdot j}$
and $\xi^-=\lim_{j\to-\infty} K_{\cdot j}$, which are
given by $\xi^+_i=i$ and $\xi^-_i=-i$. 
Thus, the Martin space $\sM$ is homeomorphic to $\Zb:=\Z\cup\{\pm\infty\}$
equipped with the usual topology.
Since both $\xi^+$ and $\xi^-$ 
are limits of almost-geodesics, $\sMin=\{\xi^+,\xi^-\}$.
Theorem~\ref{poisson-martin} says that every harmonic vector
is of the form $\alpha \xi^+ \oplus \beta\xi^- $,
for some $\alpha,\beta\in \rmax$.
\end{example}

\begin{example}\label{ex-3}
Consider $S:=\Z\times\Z$ and the operator $A$ given
by $A_{(i,j),(i,j\pm1)}=-1$ and $A_{(i,j),(i\pm1,j)}=-1$,
for each $i,j\in\Z$,
with all other entries equal to $-\infty$.
We choose the basepoint $(0,0)$.
We represent the graph of $A$
with the same conventions as in Example~\ref{ex-z}:
\begin{center}
\begin{picture}(0,0)%
\includegraphics{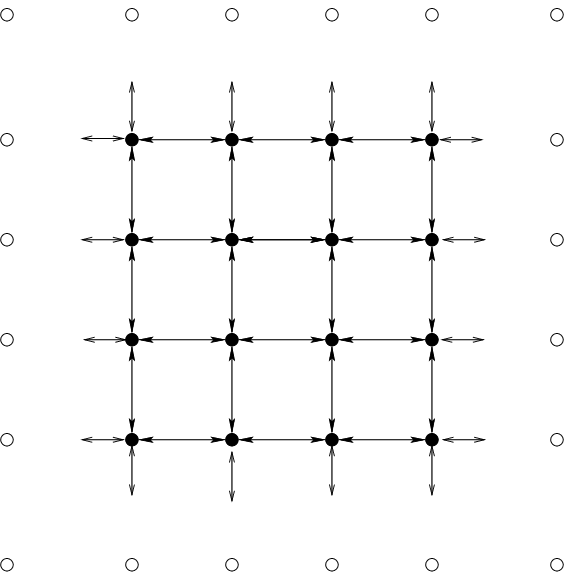}%
\end{picture}%
\setlength{\unitlength}{1579sp}%
\begingroup\makeatletter\ifx\SetFigFont\undefined%
\gdef\SetFigFont#1#2#3#4#5{%
  \reset@font\fontsize{#1}{#2pt}%
  \fontfamily{#3}\fontseries{#4}\fontshape{#5}%
  \selectfont}%
\fi\endgroup%
\begin{picture}(6766,7352)(2318,-9137)
\end{picture}%

\end{center}
For all $i,j,k,l\in\Z$,
\[
A^*_{(i,j),(k,l)} = - |i-k| - |j-l| \enspace .
\]
Note that this is the negative of the
distance in the $\ell_1$ norm between $(i,j)$ and $(k,l)$.
The matrix $A$ can be decomposed as $A=A_1\otimes I\oplus I\otimes A_2$,
where $A_1,A_2$ are two copies of the matrix of Example~\ref{ex-z},
and $I$ denotes the $\Z\times \Z$ identity matrix (recall that $\otimes$
denotes the tensor product of matrices, 
see Section~\ref{sec-decomp} for details).
The vector $\pi$ can be written as $\pi_1\otimes \pi_2$,
with $\pi_1=(A_1)^*_{0,\cdot}$ and 
$\pi_2=(A_2)^*_{0,\cdot}$. Hence, Proposition~\ref{prop-tensor}
shows that the Martin space of $A$ is homeomorphic to 
the Cartesian product of two copies of the Martin space of Example~\ref{ex-z},
in other words, that there is an homeomorphism from $\sM$ to $\Zb\times \Zb$.
Proposition~\ref{prop-tensor} also shows that
the same homeomorphism sends
$\sK$ to $\Z\times\Z$ and the minimal Martin space
to $(\{\pm\infty\}\times \Zb)
\,\cup\, (\Zb\times \{\pm\infty\})$. 
Thus, the Martin boundary and the minimal Martin space are the same.
This example may be considered to be the max-plus analogue of the
random walk on the $2$-dimensional integer lattice.
The Martin boundary for the latter (with respect to eigenvalues strictly
greater than the spectral radius) is known~\cite{neyspitzer}
to be the circle. 
\end{example}

\begin{example}\label{ex-q}
Let $S=\mathbb{Q}$ and $A_{ij}=-|i-j|$. Choosing $0$ to be the basepoint,
we get $K_{i j}=-|i -j| +|j|$ for all $j\in \mathbb{Q}$.
The Martin boundary $\sB$ consists of the functions
$i\mapsto -|i -j|+|j|$ with $j\in \R\setminus\mathbb{Q}$, together with the 
functions $i\mapsto i$ and $i\mapsto -i$. 
The Martin space $\sM$ is homeomorphic
to $\overline{\R}:=\R\cup\{\pm\infty\}$ equipped with its usual topology.
\end{example}
\begin{example}\label{ex-triangle}
We give an example of a complete locally compact metric space $(S,d)$
such that the canonical injection from $S$ to the Martin space $\sM$
is not an embedding,  and
such that the Martin boundary $\sB=\sM\setminus\sK$ is not
closed.
Consider $S=\set{(i,j)}{i\geq j\geq 1}$
and the operator $A$ given by 
\[
A_{(i,j),(i+1,j)}=A_{(i+1,j),(i,j)}=-1, \mrm{ for } i\geq j\geq 1,
\]
\[
A_{(i,j),(i,j+1)}=A_{(i,j+1),(i,j)}=-2, \mrm{ for }  i-1\geq j\geq 1,
\]
\[
A_{(1,1),(i,i)}=A_{(i,i),(0,0)}=-1/i, \mrm{ for } i\geq 2,
\]
with all other entries equal to $-\infty$.  We choose
the basepoint $(1,1)$.
The graph of $A$ is depicted in the following diagram:
\begin{center}
\begin{picture}(0,0)%
\includegraphics{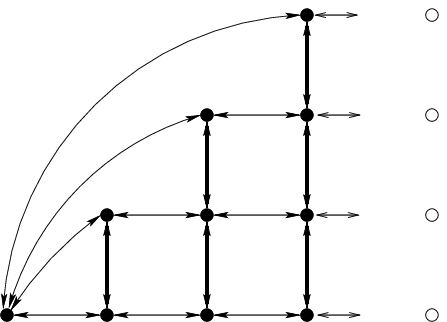}%
\end{picture}%
\setlength{\unitlength}{1579sp}%
\begingroup\makeatletter\ifx\SetFigFont\undefined%
\gdef\SetFigFont#1#2#3#4#5{%
  \reset@font\fontsize{#1}{#2pt}%
  \fontfamily{#3}\fontseries{#4}\fontshape{#5}%
  \selectfont}%
\fi\endgroup%
\begin{picture}(5266,3766)(3818,-7344)
\put(4276,-6886){\makebox(0,0)[lb]{\smash{{\SetFigFont{5}{6.0}{\rmdefault}{\mddefault}{\updefault}{\color[rgb]{0,0,0}$-1/2$}%
}}}}
\put(5251,-5536){\makebox(0,0)[lb]{\smash{{\SetFigFont{5}{6.0}{\rmdefault}{\mddefault}{\updefault}{\color[rgb]{0,0,0}$-1/3$}%
}}}}
\put(5701,-4336){\makebox(0,0)[lb]{\smash{{\SetFigFont{5}{6.0}{\rmdefault}{\mddefault}{\updefault}{\color[rgb]{0,0,0}$-1/4$}%
}}}}
\end{picture}%

\end{center}
We are using the same conventions as before.
The arcs with weight $-2$ are drawn in bold.
One can check
that 
\[
A^*_{(i,j),(k,\ell)}=
   \max\big(-|i-k|-2|j-\ell|, -(i-j)-(k-\ell)-\phi(j)-\phi(\ell) \big)
\]
where $\phi(j)=1/j$ if $j\geq 2$, and $\phi(j)=0$ if $j=1$. 
In other words, an optimal path from $(i,j)$ to $(k,\ell)$ is either
an optimal path for the metric of the weighted $\ell_1$ norm
$(i,j)\mapsto |i|+2|j|$, or a path consisting of an horizontal move
to the diagonal point $(j,j)$, followed by moves from $(j,j)$ to $(1,1)$,
from $(1,1)$ to $(\ell,\ell)$, and by an horizontal move from
$(\ell,\ell)$ to $(k,\ell)$.
Since $A$ is symmetric and $A^*$ is zero only on the diagonal,
$d((i,j),(k,\ell)):=-A^*_{(i,j),(k,\ell)}$ is a metric on $S$.
The metric space $(S,d)$ is complete since any Cauchy sequence is either
ultimately constant or converges to the point $(1,1)$.
It is also locally compact since any point distinct from
$(1,1)$ is isolated, whereas the point
$(1,1)$ has the basis of neighbourhoods consisting of the compact sets
$V_j=\set{(i,i)}{i\geq j}\cup\{(1,1)\}$, for $j\geq 2$.

If $((i_m,j_m))_{m\geq 1}$ is any sequence of elements of $S$
such that both $i_m$ and $j_m$ tend to infinity, then,
for any $(k,\ell)\in S$,
\[
A^*_{(k,\ell),(i_m,j_m)}=A^*_{(k,\ell),(1,1)}
A^*_{(1,1),(i_m,j_m)}
\qquad\text{for $m$ large enough.}
\]
(Intuitively, this is related to the fact that,
for $m$ large enough,
every optimal path from $(k,\ell)$ to $(i_m,j_m)$
passes through the point $(1,1)$).
It follows that $K_{\cdot,(i_m,j_m)}$ converges to $K_{\cdot,(1,1)}$
as $m\to\infty$.
However, the sequence $(i_m,j_m)$ does not converge to the point $(1,1)$
in the metric topology unless $i_m=j_m$ for $m$ large enough.
This shows that the map $(i,j)\to K_{\cdot,(i,j)}$
is not an homeomorphism from $S$ to its image.

The Martin boundary consists of the points
$\xi^1,\xi^2,\ldots$, obtained as limits of horizontal half-lines, which
are almost-geodesics. We have
\[
\xi^{\ell}_{(i,j)}:=\lim_{k\to \infty} K_{(i,j),(k,\ell)}=
\max\big(
i-\ell -2|j-\ell| + \phi(\ell) ,
-(i-j)-\phi(j)
\big) \enspace .
\]
The functions $\xi^\ell$ are all distinct because $i\mapsto \xi^\ell_{(i,i)}$
has a unique maximum attained at $i=\ell$. The functions $\xi^\ell$
do not belong to $\sK$ because $\xi^{\ell}_{(3j,j)}
=j+\ell+\phi(\ell)\sim j$ as $j$ tends to infinity, whereas for any $w\in \sK$,
$w_{(3j,j)}=-2j-\phi(j)\sim -2j$ as $j$ tends to infinity,.
The sequence $\xi^\ell$ converges to $K_{\cdot,(1,1)}$ as $\ell$ tends
to infinity, which shows
that the Martin boundary $\sB=\sM\setminus \sK$ is not closed.
\end{example}
\begin{example}\label{ex-corm-1}
We next give an example of a Martin space having
a boundary point which is not an extremal generator.
The same example
has been found independently by
Webster and Winchester~\cite{webster}.
Consider $S:=\mathbb{N}\times\{0,1,2\}$ and the operator $A$ given by 
\begin{equation*}
A_{(i,j),(i+1,j)}=A_{(i+1,j),(i,j)}=A_{(i,1),(i,j)}=A_{(i,j),(i,1)}=-1,
\end{equation*}
for all $i\in\N$ and $j\in\{0,2\}$,
with all other entries equal to $-\infty$.
We choose $(0,1)$ as basepoint, so that 
$\pi:=A^*_{(0,1),\cdot}$ is such that
$\pi_{(i,j)}=-(i+1)$ if $j=0$ or $2$, and 
$\pi_{(i,j)}=-(i+2)$ if $j=1$ and $i\neq 0$.
The graph associated to the matrix $A$ is depicted in
the following diagram, with the same conventions
as in the previous example.
\begin{center}
\begin{picture}(0,0)%
\includegraphics{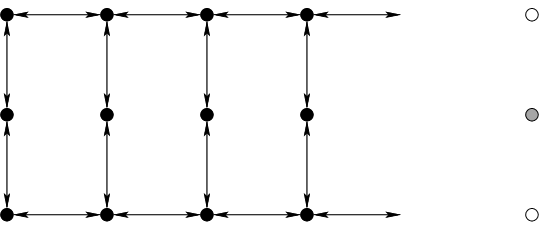}%
\end{picture}%
\setlength{\unitlength}{1579sp}%
\begingroup\makeatletter\ifx\SetFigFont\undefined%
\gdef\SetFigFont#1#2#3#4#5{%
  \reset@font\fontsize{#1}{#2pt}%
  \fontfamily{#3}\fontseries{#4}\fontshape{#5}%
  \selectfont}%
\fi\endgroup%
\begin{picture}(6466,3227)(2618,-6437)
\end{picture}%

\end{center}
There are three boundary points.
They may be obtained by taking the limits
\begin{equation*}
\xi^0:=\lim_{i\to\infty}K_{\cdot,(i,0)}, \qquad
\xi^1:=\lim_{i\to\infty}K_{\cdot,(i,1)}, \qquad
\mbox{and }
\xi^2:=\lim_{i\to\infty}K_{\cdot,(i,2)}.
\end{equation*}
Calculating, we find that
\begin{equation*}
\xi^0_{(i,j)}=i-j+1, \qquad
\xi^2_{(i,j)}=i+j-1, \qquad
\mbox{and }
\xi^1= \xi^0 \oplus \xi^2.
\end{equation*}
We have $H(\xi^0,\xi^0)=H(\xi^2,\xi^2)=H(\xi^2,\xi^1)=H(\xi^0,\xi^1)=0$.
For all other pairs $(\xi',\xi)\in\sB\times\sB$, we have $H(\xi',\xi)=-2$.
Therefore, the minimal Martin boundary is $\sMin=\{\xi^0,\xi^2\}$,
and there is a non-extremal boundary point, $\xi^1$, represented
above by a gray circle.
The sequences $((i,0))_{i\in \N}$ and $((i,2))_{i\in \N}$
are almost-geodesics,
while it should be clear from the diagram
that there are no almost-geodesics converging to $\xi^1$.
So this example provides an illustration of Propositions~\ref{semigeo-conv}
and~\ref{lemma-mr-geo}. 
\end{example}

\begin{example}\label{ex-noncompact}
Finally, we will give an example of a non-compact minimal Martin space.
Consider $S:=\N\times\N\times\{0,1\}$ and the operator $A$ given
by
\begin{align*}
 A_{(i,j,k),(i,j+1,k)} =A_{(i,j+1,k),(i,j,k)}&=-1,
\quad\mrm{for all }i,j\in\N \mrm{ and }k\in\{0,1\}, \\
A_{(i,j,k),(i,j,1-k)} &= -1,
\quad\mbox{for all } i\in\N, j\in\N\setminus\{0\}\mrm{ and }k\in\{0,1\}
, \\
A_{(i,0,k),(i,0,1-k)} &= -2,
\quad\mbox{for all } i\in\N\mrm{ and } k\in\{0,1\}, \\
A_{(i,0,k),(i+1,0,k)} =A_{(i+1,0,k),(i,0,k)}&=-1,
\quad\mbox{for all } i\in\N\mrm{ and } k\in\{0,1\}, \\
\end{align*}
with all other entries equal to $-\infty$.
We take $\pi:=A^*_{(0,0,0),\cdot}$.
With the same conventions as in Examples~\ref{ex-3} and~\ref{ex-corm-1},
the graph of $A$ is
\begin{center}
\begin{picture}(0,0)%
\includegraphics{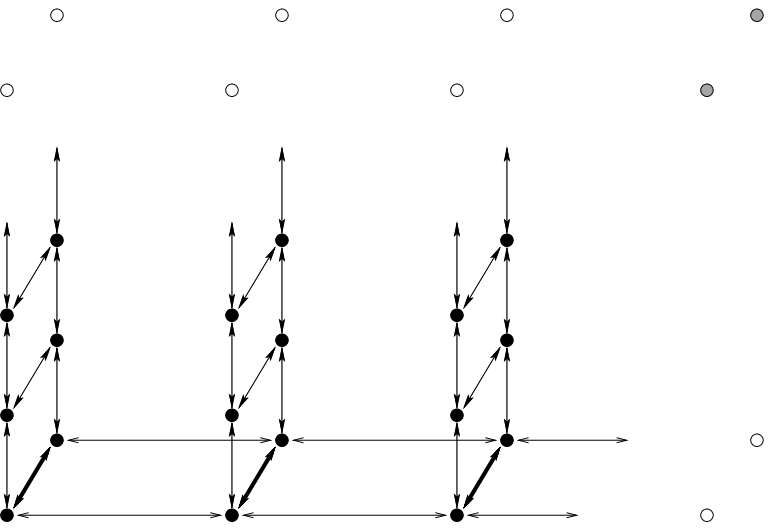}%
\end{picture}%
\setlength{\unitlength}{1579sp}%
\begingroup\makeatletter\ifx\SetFigFont\undefined%
\gdef\SetFigFont#1#2#3#4#5{%
  \reset@font\fontsize{#1}{#2pt}%
  \fontfamily{#3}\fontseries{#4}\fontshape{#5}%
  \selectfont}%
\fi\endgroup%
\begin{picture}(9166,6902)(1118,-10487)
\end{picture}%

\end{center}
Recall that arcs of weight $-1$ are drawn with thin lines whereas those
of weight $-2$ are drawn in bold.

For all $(i,j,k),(i',j',k')\in S$,
\[
A^*_{(i,j,k),(i',j',k')} = - |k'-k| - |i'-i| -
|j'-j| \ind{i=i'} -(j+j') \ind{i\neq i'}-\ind{j=j'=0,\, k\neq k'} \enspace ,
\]
where $\ind{E}$ takes the value $1$ when condition $E$ holds,
and $0$ otherwise.
Hence, 
\begin{align*}
{K}_{(i,j,k),(i',j',k')} =& k'- |k'-k| +i'- |i'-i| +
j'- |j'-j| \ind{i=i'} -(j+j') \ind{i\neq i'} \nonumber \\
&+\ind{j'=0,k'=1}-\ind{j=j'=0,\, k\neq k'} \enspace .
\end{align*}
By computing the limits of $K_{\cdot,(i',j',k')}$
when $i'$ and/or $j'$ go to $+\infty$,
we readily check that the Martin boundary is composed of
the vectors 
\begin{align*}
\xi^{i',\infty,k'} &:=\lim_{j'\to\infty}K_{\cdot,(i',j',k')}, \\
\xi^{\infty,\infty,k'} &:=\lim_{i',j'\to\infty} K_{\cdot,(i',j',k')} \\
\xi^{\infty,0,k'} &:=\lim_{i'\to\infty}K_{\cdot,(i',0,k')}.
\end{align*}
where the limit in $i$ and $j'$ in the second line can be taken in either order.
Note that $\lim_{i'\to\infty}K_{\cdot,(i',j',k')}=\xi^{\infty,\infty,k'}$ 
for any $j'\in  \N\setminus\{0\}$ and $k'\in\{{0,1}\}$.
The minimal Martin space is composed of the vectors
$\xi^{i',\infty,k'}$ and $\xi^{\infty,0,k'}$ with $i'\in\N$ and $k'\in\{0,1\}$.
The two boundary points $\xi^{\infty,\infty,0}$ and $\xi^{\infty,\infty,1}$ 
are non-extremal and have representations
\begin{align*}
\xi^{\infty,\infty,0} &=\xi^{\infty,0,0}\oplus-3 \xi^{\infty,0,1} \enspace,
\\
\xi^{\infty,\infty,1} &=\xi^{\infty,0,0}\oplus-1 \xi^{\infty,0,1}.
\end{align*}
For $k'\in\{0,1\}$,  the sequence
$(\xi^{i',\infty,k'})_{i\in\N}$ converges
to $\xi^{\infty,\infty,k'}$ as $i$ goes to infinity.
Since this point is not in $\sMin$, we see that $\sMin$ is not compact.
\end{example}
\section{Tightness and existence of harmonic vectors}\label{sec-tight}
\label{subsec-exists}
We now show how the Martin boundary can be used 
to obtain existence results for eigenvectors.
As in~\cite{AGW-s}, we restrict our attention
to the case where $S$ is equipped with the discrete topology.
We say that a vector $u\in \rmax^S$ is
\new{$A$-tight} if, for all $i\in S$ and $\beta\in\R$,
the super-level set $\set{j\in S}{A_{ij}u_j\geq \beta}$
is finite. 
We say that a family of vectors $\{u^\ell\}_{\ell\in L}\subset
\rmax^S$ is $A$-tight if $\sup_{\ell\in L} u^\ell$
is $A$-tight. The notion of tightness is motivated by the
following property.
\begin{lemma}\label{lem-fond}
If a net $\{u^\ell\}_{\ell\in L}\subset \rmax^S$
is $A$-tight and converges pointwise to $u$, then $Au^\ell$
converges pointwise to $Au$.
\end{lemma}
\proof
This may be checked elementarily, 
or obtained as a special case of general results
for idempotent measures~\cite{12,aqv95,DENSITE,tolya01}
or, even more generally, capacities~\cite{brienv}.
We may regard $u$ and $u^l$ as the densities of the idempotent measures
defined by
\begin{equation*}
Q_u(J)=\sup_{j\in J} u_j
\qquad\mbox{and}\qquad
Q_{u^l}(J)=\sup_{j\in J} u^l_j
\enspace  ,
\end{equation*}
for any $J\subset S$.
When $S$ is equipped with the discrete topology, pointwise convergence
of $(u^\ell)_{\ell\in L}$ is equivalent to convergence
in the hypograph sense of convex analysis.
It is shown in~\cite{aqv95} that this is then equivalent to
convergence of $(Q_{u^l})_{\ell\in L}$
in a sense analogous to the vague convergence of
probability theory. It is also shown that, when combined with the tightness
of $(u_l)_{\ell\in L}$,
this implies convergence in a sense analogous to weak convergence.
The result follows as a special case.
\qed
\begin{prop}\label{prop-cs}
Assume that $S$ is infinite and that 
the vector $\pi^{-1}:=(\pi_i^{-1})_{i\in S}$
is $A$-tight. Then, some element of $\sM$ is harmonic
and, if $\zero\not\in\sM$, then $\sMin$ is non-empty.
Furthermore, each element of $\sB$ is harmonic. 
\end{prop}
\begin{proof}
Since $S$ is infinite, there exists an injective map $n\in \N\mapsto
i_n\in S$. Consider the sequence $(i_n)_{n\in \N}$.
Since $\sM$ is compact, it has a subnet 
$(j_k)_{k\in D}$, $j_k:=i_{n_k}$ such that 
$\{K_{\cdot j_k}\}_{k\in K}$ converges to some $w\in \sM$.
Let $i\in S$. Since $(AA^*)_{ij}=A^+_{ij}=A^*_{ij}$ for all
$j\neq i$, we have 
\begin{align*}
(AK_{\cdot j_k})_i= K_{i j_k}
\end{align*}
when $j_k \neq i$.
But, by construction, the net $(j_k)_{k\in D}$ is eventually
in $S\backslash\{i\}$ and so we may pass to the limit,
obtaining $\lim_{k\in K} A K_{\cdot j_k}=w$.
Since $\pi^{-1}$ is $A$-tight, it follows
from~\eqref{ine-K} that the family
$(K_{\cdot j})_{j\in S}$
is $A$-tight. Therefore, by Lemma~\ref{lem-fond},
we get $w= Aw$.
If $\zero\not\in\sM$, then $\sH$ contains
a non-zero vector, and applying the representation
formula~\eqref{Hequal1} to this vector, we see that $\sMin$ cannot be empty.

It remains to show that $\sB\subset\sH$.
Any $w\in\sB$ is the limit of a net $\{K_{\cdot j_k}\}_{k\in D}$.
Let $i\in S$.
Since $w\neq K_{\cdot i}$, the net $\{K_{\cdot j_k}\}_{k\in D}$
is eventually in some neighbourhood of $w$ not containing
$K_{\cdot i}$.
We deduce as before that $w$ is harmonic.
\end{proof}
\begin{corollary}[Existence of harmonic vectors]
\label{cor-a-tight}
Assume that $S$ is infinite, 
that $\pi=A^*_{b\cdot}\in\R^S$ for some $b\in S$,
and that $\pi^{-1}$ is $A$-tight. Then, $\sH$ contains
a non-zero vector.
\end{corollary}
\begin{proof}
We have $K_{bj}=\unit$ for all $j\in S$ and hence,
by continuity, $w_b=\unit$ for all $w\in\sM$.
In particular,  $\sM$ does not contain $\zero$.
The result follows from an application of the proposition.
\end{proof}
We finally derive a characterisation of the spectrum of $A$.
We say that $\lambda$ is a {\em (right)-eigenvalue} of $A$ if
$Au=\lambda u$ for some vector $u$ such that $u\neq \zero$.
\begin{corollary}\label{cor-a-irred}
Assume that $S$ is infinite,
$A$ is irreducible, and for each $i\in S$, there are only finitely
many $j\in S$ with $A_{ij}>\zero$.
Then the set of right eigenvalues of $A$ is $[\rho(A),\infty[$.
\end{corollary}
\begin{proof}
Since $A$ is irreducible, no eigenvector of $A$ can have a component
equal to $\zero$.
It follows from~\cite[Prop.~3.5]{dudnikov} that
every eigenvalue of $A$ must be greater than or equal to~$\rho(A)$.

Conversely, for all $\lambda\geq \rho(A)$, we have
$\rho(\lambda^{-1} A)\leq \unit$.
Combined with the irreducibility of $A$,
this implies~\cite[Proposition 2.3]{AGW-s} that
all the entries of $(\lambda^{-1} A)^*$ are finite.
In particular, for any $b\in S$, the vector
$\pi:=(\lambda^{-1} A)^*_{b \cdot}$ is in $\R^S$.
The last of our three assumptions ensures that
$\pi^{-1}$ is $(\lambda^{-1} A)$-tight
and so, by Corollary~\ref{cor-a-tight},
$(\lambda^{-1} A)$ has a non-zero harmonic vector.
This vector will necessarily be an eigenvector of $A$ with
eigenvalue $\lambda$.
\end{proof}
\begin{example}
The following example shows that when $\pi^{-1}$ is not $A$-tight,
a Martin boundary point need not be an eigenvector.
Consider $S:=\N$ and the operator $A$ given by 
\begin{align*}
A_{i,i+1}=A_{i+1,i}:=-1
\qquad\mbox{and}\qquad
A_{0i}&:=0
\qquad \mbox{for all $i\in \N$},
\end{align*}
with all other entries of equal to $-\infty$.
We take $\pi:=A^*_{0,\cdot}$.
With the same conventions as in Example~\ref{ex-corm-1},
the graph of $A$ is
\begin{center}
\begin{picture}(0,0)%
\includegraphics{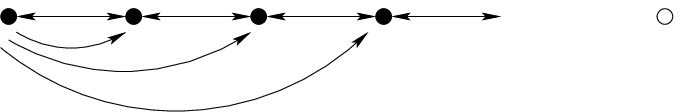}%
\end{picture}%
\setlength{\unitlength}{1973sp}%
\begingroup\makeatletter\ifx\SetFigFont\undefined%
\gdef\SetFigFont#1#2#3#4#5{%
  \reset@font\fontsize{#1}{#2pt}%
  \fontfamily{#3}\fontseries{#4}\fontshape{#5}%
  \selectfont}%
\fi\endgroup%
\begin{picture}(6466,1217)(2618,-4577)
\put(4201,-4036){\makebox(0,0)[lb]{\smash{{\SetFigFont{6}{7.2}{\rmdefault}{\mddefault}{\itdefault}{\color[rgb]{0,0,0}$0$}%
}}}}
\put(5401,-4186){\makebox(0,0)[lb]{\smash{{\SetFigFont{6}{7.2}{\rmdefault}{\mddefault}{\itdefault}{\color[rgb]{0,0,0}$0$}%
}}}}
\put(3151,-3861){\makebox(0,0)[lb]{\smash{{\SetFigFont{6}{7.2}{\rmdefault}{\mddefault}{\itdefault}{\color[rgb]{0,0,0}$0$}%
}}}}
\end{picture}%

\end{center}
We have $A^*_{i,j}=\max(-i,-|i-j|)$
and $\pi_i=0$ for all $i,j\in \N$.
There is only one boundary point, 
$b :=\lim_{k\to\infty} K_{\cdot k}$,
which is given by $b_i=-i$ for all $i\in\N$. 
One readily checks that $b$ is not an harmonic vector
and, in fact, $A$ has no non-zero harmonic vectors.
\end{example}
\section{Eigenvectors of Lax-Oleinik semigroups and Busemann points of normed spaces}\label{sec-lax}
We now use the Martin boundary to solve
a class of continuous-time deterministic optimal control problems.
Consider the value function $v$ defined  by:
\[
v(t,x):= \sup_{X(\cdot),\; X(0)=x} 
\phi (X(t)) -\int_0^t L(\dot X(s)) \dd s   \enspace .
\]
Here, $x$ is a point in $\R^n$, $t$ is a nonnegative real number,
the {\em Lagrangian} $L$ is a Borel measurable map $\R^n \to \R\cup\{+\infty\}$, bounded from below, the {\em terminal reward}
$\phi$ is an arbitrary map $\R^n\to \R\cup\{-\infty\}$,
and the supremum is taken over all absolutely continuous
functions $X: [0,t]\to \R^n$ such that $X(0)=x$. 
This is a special case of the classical Lagrange problem
of calculus of variations.

The {\em Lax-Oleinik semigroup} $(T^t)_{t\geq 0}$ is 
composed of the maps $T^t$ sending
the value function at time $0$, $v(0,\cdot)=\phi$ to the value 
function at time $t$, $v(t,\cdot)$. 
The semigroup property $T^{t+s}=T^t\circ T^s$
follows from the dynamic programming principle.
The kernel of the operator $T^t$
is given by 
\[
(x,y)\mapsto 
T^t_{x,y}= \sup_{X(\cdot),\; X(0)=x,\;X(t)=y} -\int_0^t L(\dot X(s)) \dd s 
\enspace,
\]
where the supremum is taken over all absolutely continuous
functions $X:[0,t]\to \R^n$ such that $X(0)=x$ and $X(t)=y$. 

The classical Hopf-Lax formula states that
\[
T^t_{x,y} = -t\operatorname{co} L\Big(\frac{y-x}{t}\Big) ,\qquad \text{for }t>0\enspace,
\]
where $\operatorname{co} L$ denotes the convex lower semicontinuous hull of $L$.
This is proved, for instance, in~\cite[\S 3.3, Th.~4]{evans} when $L$ is convex
and finite valued, and when the curves $X(\cdot)$ are required to be continuously differentiable. The extension to the present setting is not difficult.

Since $T^t$ only depends on $\operatorname{co} L$, we shall assume
that $L$ is convex, lower semicontinuous,
and bounded from below. Moreover, we shall always
assume that $L(0)$ is finite.

We say that a function $u:\R^n\to \R\cup\{-\infty\}$,
not identically $-\infty$, is an {\em eigenvector} 
of the semigroup $(T^t)_{t\geq 0}$ with
{\em eigenvalue} $\lambda$ if 
\[
T^t u=u+\lambda t,\;\text{ for all }t> 0 \enspace .
\]
We shall say that $u$ is {\em extremal} if it is an extremal
generator of the eigenspace of the semigroup $(T^t)_{t\geq 0}$
with eigenvalue $\lambda$,
meaning that $u$ cannot be written as the supremum of two eigenvectors
with the same eigenvalue that are both different from it.

One easily checks, using the convexity of $L$, that for all $t>0$,
the maximal circuit mean of the operator $T^t$ is given by
\[
\rho(T^t) = -tL(0) \enspace.
\]
By Proposition~3.5 of~\cite{dudnikov} or Lemma~2.2 of~\cite{AGW-s},
any eigenvalue $\mu$ of $T^t$ must satisfy $\mu\geq \rho(T^t)$,
and so any eigenvalue $\lambda$ of the semigroup $(T^t)_{t\geq 0}$
satisfies
\[ \lambda \geq -L(0)\enspace .\]
We denote by $\zeta(x)$ the one sided directional derivative of $L$ at
the origin in the direction $x$:
\begin{align}\label{e-dirder}
\zeta(x)=\lim_{t\to 0^+} t^{-1}(L(tx)-L(0))=\inf_{t>0}
t^{-1}(L(tx)-L(0))\in\R\cup\{\pm\infty\} \enspace,
\end{align}
which always exists since $L$ is convex.
\begin{prop}\label{cor-lax1}
Assume that $\zeta$ does not take the value $-\infty$.
Then, the eigenvectors of the Lax-Oleinik semigroup $(T^t)_{t\geq 0}$
with eigenvalue $-L(0)$ are 
precisely the functions $u:\R^n\to \R\cup\{-\infty\}$, not identically
$-\infty$, such that 
\begin{align}\label{e-repr-L}
-\zeta(y-x)+ u(y)\leq u(x) \enspace ,
\text{ for all }x,y\in \R^n\enspace . 
\end{align}
Moreover, when $\zeta$ only takes finite values,
the extremal eigenvectors with eigenvalue $-L(0)$
are of the form $c+w$, where
$c \in \R$ and $w$ belongs to the minimal Martin space of
the kernel $(x,y)\mapsto -\zeta(y-x)$ with respect to any basepoint.
\end{prop}
\begin{proof}
Let us introduce the kernels
\[
A_s:=T^s+sL(0) ,\text{ for all }s\geq 0.
\]
Using the Hopf-Lax formula, we get
\begin{align*}
(A_s)^+_{xy}&=\sup_{k\in \N\setminus\{0\}} -ks L\Big(\frac{y-x}{ks}\Big)+ks L(0)
\enspace .
\end{align*}
Using~\eqref{e-dirder}
and the fact that $\zeta(0)=0$, we deduce that
\begin{align}\label{e-*}
(A_s)^*_{xy} =(A_s)^+_{xy}=-\zeta(y-x) \enspace.
\end{align}
The eigenvectors of the semigroup $(T^t)_{t\geq 0}$ are
precisely the functions that are harmonic with respect
to all the kernels $A_s$, with $s>0$. 
Since $(A_s)_{xx}=0$ for all $x\in\R^n$, the harmonic
and super-harmonic functions of $A_s$ coincide.
It follows from Proposition~\ref{superharm} that
$u$ is a super-harmonic function of $A_s$ if and only if
$u\geq A_s^* u$. Since the latter condition can be written
as~\eqref{e-repr-L} and is independent of $s$,
the first assertion of the corollary is proved.

By~\eqref{e-*}, when $\zeta$ is finite, any point can
be taken as the basepoint.
The kernels $A_s$ and $(x,y)\mapsto -\zeta(y-x)$
have the same Martin and minimal Martin spaces
with respect to any given basepoint,
and so the final assertion of the corollary follows from
Theorem~\ref{th-mr-ext}.
\end{proof}
\begin{remark}
When $\partial L(0)$, the subdifferential of $L$
at the origin, is non-empty, $\zeta$
does not take the value $-\infty$.
This is the case when the origin is in the relative interior
of the domain of $L$. Then, $\zeta$ coincides with the support function of
$\partial L(0)$:
\[
\zeta(x)= \sup_{y\in \partial L(0)}{y\cdot x}  ,\qquad
\text{for all } x\in \R^n \enspace ,
\]
see~\cite[Th.~23.4]{ROCK}.
If in addition the origin is in the interior of the
domain of $L$, then $\partial L(0)$ is non-empty and compact,
and so the function $\zeta$ is everywhere finite.
\end{remark}
\begin{corollary}\label{cor-norm1}
When $\zeta$ is a norm on $\R^n$, the extremal eigenvectors
with eigenvalue $-L(0)$ of the Lax-Oleinik semigroup $(T^t)_{t\geq 0}$
are precisely the functions
$x\mapsto c-\zeta(y-x)$, where
$c\in\R$ and $y\in \R^n$, together with the functions $c+w$,
where $c\in \R$ and $w$ is a Busemann point 
of the normed space $(\R^n,\zeta)$.
\end{corollary}
\begin{proof}
This follows from Proposition~\ref{cor-lax1} and Corollary~\ref{cor-rieffel}.
\end{proof}
\begin{remark}
The map $\zeta$ is a norm when the origin
is in the interior of the domain of $L$
and the subdifferential $\partial L(0)$ is symmetric, meaning that $p\in \partial L(0)$ implies $-p\in \partial L(0)$. When $\zeta$ is a norm, 
condition~\eqref{e-repr-L}
means that $u$ is Lipschitz-continuous with respect
to $\zeta$ or that $u$ is identically $-\infty$.
\end{remark}
We next study the eigenspace of $(T^t)_{t\geq 0}$ for an eigenvalue
$\lambda>-L(0)$ in the special case where $L$ is of the form
\begin{align*}
L(x)= \frac{\|x\|}p^p \enspace ,
\end{align*}
where $\|\cdot\|$ is an arbitrary norm on $\R^n$ and $p> 1$.
For all $\lambda>0$, we set
\[
\vartheta_\lambda :=(q\lambda)^{\frac 1 q}
\quad\text{ where }\frac 1p+\frac 1q=1\enspace .
\]
\begin{theorem}\label{cor-final}
Let $s>0$ and $\lambda>0$. Any eigenvector of 
$T^s$ with eigenvalue $\lambda s$ is an
eigenvector of the Lax-Oleinik
semigroup $(T^t)_{t\geq 0}$ with eigenvalue
$\lambda$. Such an eigenvector can be written as
\begin{align}\label{e-rep-spec}
u= \sup_{w\in \mbuu} \nu(w) +  \vartheta_\lambda w \enspace,
\end{align}
where $\mbuu$ denotes the set of Busemann points
of the normed space $(\R^n,\|\cdot\|)$
and $\nu$ is an arbitrary map $\mbuu\to \R\cup\{-\infty\}$
bounded from above. The maximal map $\nu$ satisfying~\eqref{e-rep-spec}
is given by $\mu_u$.
Moreover, the extremal eigenvectors with eigenvalue
$\lambda$ are of the form $c+\vartheta_\lambda w$, where $c\in\R$ and $w\in\mbuu$.
\end{theorem}
This theorem follows from Theorem~\ref{poisson-martin},
Theorem~\ref{th-mr-ext2}, and the next lemma.
\begin{lemma}\label{lemma-technical}
For all $s>0$, the minimal Martin space of the kernel 
$A_s:=T^s-s\lambda$, with respect to any basepoint, coincides
with the set of functions $\vartheta_\lambda w$, where $w$ is
a Busemann point of the normed space $(\R^n,\|\cdot\|)$
equipped with the same basepoint.
\end{lemma}
\begin{proof}
For all $x,y\in\R^n$, 
we set
\[ \psi(t):= -t^{1-p}L(y-x) - t \lambda \enspace .
\]
It follows from the Hopf-Lax formula that
\begin{align}
(A_s)^+_{xy}= \sup_{k\in \N\setminus\{0\}}
\psi(ks) \enspace.\label{e-opt}
\end{align}
Since $\psi$ is concave, the supremum of $\psi(t)$ over all $t>0$ is attained
at the point $\bar t$ such that
\begin{align*}
\psi'(\bar t)=\bar t^{-p}(p-1)L(y-x)- \lambda =0 \enspace.
\end{align*}
It follows that
\[
\psi(\bar t) = -\vartheta_{\lambda}\|y-x\| \enspace .
\]
Since $\psi$ is concave, we have $\psi(t)\geq \psi(\bar t)+\psi'(t)(t-\bar t)$,
and so, for $t\geq \bar t$, 
\begin{align*}
\psi(t)-\psi(\bar t)&=\psi(t)-\psi(\bar t)-\psi'(\bar t)(t-\bar t)\\
&\geq (\psi'(t)-\psi'(\bar t))(t-\bar t) \geq\psi''(\bar t)(t-\bar t)^2 
\end{align*}
since $\psi'$ is convex.
Let $k$ denote the smallest integer such that $\bar t\leq ks$, and let $t=ks$.
We deduce that
\begin{align*}
0\geq \psi(t)-\psi(\bar t)&\geq -p(p-1)L(y-x){\bar t}^{-1-p}(t-\bar t)^2
= -p\lambda {\bar t}^{-1}(t-\bar t)^2 \enspace .
\end{align*}
Since $\bar t \leq t\leq \bar t +s$, since $\bar t =(q\lambda)^{-1/p}\|y-x\|$,
and since 
\[ \psi(\bar t) \geq (A_s)^*_{xy}\geq (A_s)^+_{xy}\geq \psi(t)\enspace,\]
we get
\begin{align}\label{e-theta}
(A_{s})_{xy}^* = -\vartheta_{\lambda}\|y-x\| +\epsilon(\|y-x\|) \enspace ,
\end{align}
where $\epsilon$ is a function tending to $0$ at infinity.
Observe that the supremum in~\eqref{e-opt} is always
attained by an integer $k$ which can be bounded by 
an increasing function of $\|y-x\|$. Hence, for
all $x\in \R^n$ and every compact set $C$,
we can find an integer $N$ such that
$(A_s)^+_{xy}=\sup_{1\leq k\leq N}\psi(ks)$
for all $y\in C$. Since every $\psi(ks)$
is a continuous function of $y-x$, we deduce that
the map $y\mapsto (A_s)^+_{xy}$ is continuous.

Denote by
$K$ the Martin kernel of $A_s$ 
with respect to this basepoint and denote by $\sM$, $\sMin$,
and $\sK$,  the corresponding Martin space, minimal
Martin space,
and set of columns of the Martin kernel. Also, we denote by $H$
the kernel
constructed from $K$
as in Section~\ref{sec-minmartin}.
Define the kernel
$A': (x,y)\mapsto -\vartheta_\lambda\|y-x\|$.
We use $K',\sM',\sMinp, \sK'$ and $H'$ to denote
the corresponding objects constructed from $A'$.

We next show that $\sMin=\sMinp\setminus\sK'$.

An element $w$ of $\sMin$ is
the limit of a net $(K_{\cdot  y_d})_{d\in D}$. 
If the net $(y_d)_{d\in D}$ had a bounded subnet,
it would have a subnet
converging to some $y\in \R^d$. Then, by continuity
of the map $z\mapsto (A_s)^+_{\cdot z}$, 
the element $w$ would be proportional in the max-plus sense either
to $f:=(A_s)^*_{\cdot y}$ or to $g:=(A_s)^+_{\cdot y}$ (the first
case arises if the subnet is ultimately constant).
Both cases can be ruled out: we know from Proposition~\ref{mrsubh}
that an element of the minimal Martin space is harmonic,
but $f_y=0\neq
g_y=(A_s f)_y=-s\lambda\neq
(A_sg)_y=-2s\lambda$,
and so $f$ and $g$ are not harmonic.
This shows that $(y_d)_{d\in D}$ tends to infinity.

By~\eqref{e-theta}, we deduce that $K'_{\cdot y_d}$ tends to $w$.
Thus, any net $(y_d)_{d\in D}$ such that $K_{\cdot y_d}$ tends to $w$
is such that $y_d$ tends to infinity
and $K'_{\cdot y_d}$ tends to $w$.
We deduce that $w\in \sM'$ and $H'(w,w)\geq H(w,w)=\unit$,
and so, by~\eqref{prop-mr}, $\sMin\subset\sMinp\cup\sK'$. 

We proved that the columns of $(A_s)^*$ are not harmonic,
and so $\sMin\subset \sM\setminus \sK$. 
We claim that $\sMin\subset \sMinp\setminus\sK'$. 
Indeed, if a net $K_{\cdot y_d}$ converges to $w\in\sMin$,
we showed that $(y_d)_{d\in D}$ tends to infinity,
and that $K'_{\cdot  y_d}$ tends to $w$. But
$K'_{\cdot  y_d}$ cannot converge to an element
$K'_{\cdot  y}\in \sK'$ because the map sending an element of a
finite-dimensional normed space to its column of the Martin kernel
is an embedding (see~\cite[Ch.~II,\S 1]{ballmann}
for a more general result). So $w\not\in\sK'$.

Let us take now $w'\in\sMinp\setminus\sK'$.
Then, $w'$ is the limit
of some net $(K'_{\cdot  y'_d})_{d\in D'}$,
where $(y'_d)_{d\in D'}$ necessarily tends to infinity,
since otherwise, there would be a subnet of 
$(y'_d)_{d\in D'}$ converging to some $z\in\R^n$,
and so we would have $w'=K'_{\cdot z}\in\sK'$.
It follows from~\eqref{e-theta} that $w'$ is the limit
of $K_{\cdot y'_d}$, and hence $w'\in \sM$.
These properties also imply that $H'(w',w')\leq H(w',w')$.
Since $w'\in \sMinp$, we have $H'(w',w')=\unit$, and so $H(w',w')=\unit$,
and by~\eqref{prop-mr}, $w'\in \sMin\cup \sK$.
Observe that the map $z\mapsto w'_z$ is continuous
because it is a pointwise limit
of elements of $\sK'$, all of which are Lipschitz continuous
with constant $\vartheta_\lambda$ with respect to the norm $\|\cdot\|$.
For all $y\in \R^n$, the map $x\mapsto A^*_{xy}$ takes the value $0$
when $x=y$ and the value $(A_s)^+_{xy}\leq -s\lambda<0$ when $x\neq y$.
Thus, the elements of $\sK$ are not continuous, and so, $w'\not\in \sK$.
It follows that $w'\in \sMin\setminus \sK=\sMin$. 
We have shown that $\sMin=\sMinp\setminus\sK'$.

By Corollary~\ref{cor-rieffel}, $\sMinp\setminus\sK'$ is the set of Busemann
points of the normed space $(\R^n,\vartheta_\lambda\|\cdot\|)$.
These are precisely the functions of the form $\vartheta_\lambda w$,
where $w$ is a Busemann point of $(\R^n ,\|\cdot\|)$.
\end{proof}
\begin{remark}
Lemma~\ref{lemma-technical} identifies a special situation where
the minimal Martin space of $T^s-s\lambda$ is independent of $s$.
This seems related to the fact that the set of functions
of the form $x\mapsto a\|x\|^p$ with $a>0$ is stable by inf-convolution.
One may still obtain
a representation of the eigenvectors for more general semigroups
$(T^t)_{t\geq 0}$, but this requires adapting some of the present results
to the continuous-time setting. We shall present this elsewhere.
\end{remark}
\begin{example}
Consider the Euclidean norm on $\R^n$, $\|x\|:=(x\cdot x)^{1/2}$, 
and $L(x):=\|x\|^p/p$ with $p>1$. 
The set of Busemann points of the normed space $(\R^n,\|\cdot\| )$,
with respect to the basepoint $0$,
coincides with the set of functions 
\[ w: x\mapsto 
x\cdot y \enspace ,
\]
where $y$ is an arbitrary vector of norm $1$.
It follows from Theorem~\ref{cor-final}
that the extremal eigenvectors
with eigenvalue $\lambda>0$ 
of the Lax-Oleinik semigroup are of the form
$c+\vartheta_\lambda w$, with $c\in \R$, and that any eigenvector
with eigenvalue $\lambda$ is a supremum of maps of this form.
In particular, when $n=1$, there are two Busemann points,
$w^\pm (x)= \pm \vartheta_\lambda x$, and
any eigenvector $u$ with eigenvalue
$\lambda$ can be written as
\begin{align*}
x\mapsto \max(c^++\vartheta_\lambda x, c^--\vartheta_\lambda x) \enspace ,
\end{align*}
with $c^\pm\in \R\cup\{-\infty\}$. 
The Busemann points $w^\pm$ are the limits of the geodesics
$t\mapsto \pm t$, from $[0,\infty[$ to $\R$.
Hence, Proposition~\ref{prop-carac-specmes} allows us to determine
the maximal representing measure $\mu_u$, or equivalently, the maximal value
of the scalars $c^{\pm}$, as follows:
\[
c^\pm = \lim_{t \to \pm \infty} u(t)\mp \vartheta_\lambda t \enspace .
\]
In this special case, the representing measure is unique.
\end{example}
In order to give another example, we characterise the Busemann points
of a polyhedral norm. We call {\em proper face} of a polytope
the intersection of this polytope with a supporting half-space.
\begin{prop}\label{prop-carac-bus}
Let $\|\cdot\|$ denote a polyhedral norm on $\R^n$, so that
\[
\|x\|=\max_{i\in I}x'_i \cdot x\enspace ,
\]
where $(x'_i)_{i\in I}$ is the finite family
of the extreme
points of the dual unit ball.
The Martin boundary of the kernel
$(x,y)\mapsto -\|x-y\|$, taking the origin as the basepoint,
is precisely the set of functions of the form
\begin{align}\label{e-bus}
x\mapsto \min_{j\in J} x'_j\cdot (x-X)+\max_{j\in J} x'_j \cdot X\enspace,
\end{align}
where $X\in \R^n$ and $(x'_j)_{j\in J}$ is the set of extreme points
of a proper face of the dual unit ball.
Moreover, all the points
of the Martin boundary are Busemann points.
\end{prop}
\begin{proof}
Any point $f$ of the Martin boundary
is the limit of a sequence of functions 
\[
x\mapsto f^k(x) = \|X^k\|-\|X^k-x\| \enspace ,
\]
where $X^k\in \R^n$ and $\|X^k\|\to\infty$ when $k\to\infty$.
Consider the sequence of vectors 
\[ u^k=(x'_i\cdot X^k-\|X^k\|)_{i\in I}
\enspace .
\]
These vectors lie in $[-\infty,0]^I$,
which is compact and metrisable, and so,
we may assume,
by taking a subsequence if necessary,
that $u^k$ converges to some
vector $u\in [-\infty,0]^I$.
Since $I$ is finite, we may also assume,
again taking a subsequence if necessary, 
that there exists
an index $j_0\in I$ such that $x'_{j_0}\cdot X^k=\|X^k\|$
for all $k$. Let $J:=\set{i\in I}{u_i>-\infty}$. 
Observe that $J$ is non-empty since $u_{j_0}=0$. 
We have
\begin{align*}
f(x)&=\lim_{k\to\infty} f^k(x) = \lim_{k\to\infty} -\max_{i\in I} (x'_i\cdot X^k-\|X^k\|-x'_i\cdot x)\\
&=-\max_{j\in J} (u_j -x'_{j}\cdot x) \enspace .
\end{align*}
Observe that the set $E:=\set{((x'_j-x'_{j_0})\cdot X)_{j\in J}}{X\in\R^n}$
is closed, since it is a finite-dimensional vector space.
Since the vector $(u^k)_{j\in J}$ belongs
to $E$ and has a finite limit when $k\to\infty$, this
limit belongs to $E$, and so there exists some $X\in \R^n$ such that
$u_j=x'_j\cdot X-x'_{j_0}\cdot X$ for all $j\in J$.
Thus, 
\begin{align*}
f(x)&=-\max_{j\in J} x'_j\cdot (X-x) + x'_{j_0}\cdot X \enspace .
\end{align*}
Since $f(0)=0$, we have $\max_{j\in J}x'_j\cdot X=x'_{j_0}\cdot X$,
and so
\begin{align*}
f(x)&=-\max_{j\in J} x'_j\cdot (X-x) + \max_{j\in J}x'_{j}\cdot X 
\enspace,
\end{align*}
which is of the form~\eqref{e-bus}.

We now have to show that $(x'_j)_{j\in J}$ is the set of extreme points
of a face of the dual unit ball.
Let $E'$ denote the set of vectors $x'\in\R^n$ such that $x'\cdot X^k-\|X^k\|$
remains bounded when $k$ tends to infinity. This is
an affine space. Let $B'$ denote the dual unit ball.
We claim that $F':=E'\cap B'$ is an extreme subset of $B'$,
meaning that
\begin{align}
\alpha x'+(1-\alpha) y' \in F'\implies x',y'\in F',
\qquad
\text{for all }x',y'\in B'
\text{ and }0< \alpha< 1.
\label{e-extreme}
\end{align}
Indeed, let $x',y'\in B'$ and $0<\alpha<1$.
Since $x'\in B'$,
we have $x'\cdot X\leq \|X\|$ for all $X\in \R^n$.
In particular, $x'\cdot X^k-\|X^k\|\leq 0$ for all $k$.
Similarly, $y'\cdot X^k-\|X^k\|\leq 0$ for all $k$.
Since 
\begin{align*}
(\alpha x'+(1-\alpha) y')\cdot X^k -\|X^k\|
&= \alpha (x'\cdot X^k-\|X^k\|)
+(1-\alpha) (y'\cdot X^k -\|X^k\|)\\
&\leq \alpha (x'\cdot X^k -\|X^k\|)\\
&\leq 0\enspace,
\end{align*}
we deduce that $x'\cdot X^k -\|X^k\|$ is bounded
if $\alpha x'+(1-\alpha) y'\in F'$.
Similarly, $y'\cdot X^k -\|X^k\|$ is bounded.
This shows~\eqref{e-extreme}.

Let $z$ denote any accumulation point of the sequence
$\|X^k\|^{-1}X^k$. We have $F'\subset \set{x'\in B'}{x'\cdot z =1}$,
and so, $F'\neq B'$.

Since the dual ball $B'$ is a polytope,
the convex extreme subset $F'\neq B'$
is a proper face of $B'$. Therefore, the vectors $x'_i$,
with $i\in I$,
such that $x'_i\cdot X^k-\|X^k\|$ remains bounded
are precisely the $x'_i$ that belong to the proper face $F'$.
Hence, these $x'_i$ are the extreme points of the proper face $F'$.

Every proper face $F'$ of the dual ball is the intersection of the dual ball
with a supporting hyperplane, so $F'=\set{x'\in B'}{x'\cdot y=1}$ for
some $y\in B$. Observe that the set $J$ of $x_i'$ such that $x_i'\cdot y=1$
is precisely the set of extreme points of $F'$.
Consider now $X\in \R^n$ and the ray $t\mapsto X+ t y$, which is a geodesic, and a fortiori an almost-geodesic.
One readily checks that the function
$x\mapsto \|X+t y\|-\|X+t y -x\|$
converges to the function~\eqref{e-bus} when $t$ tends to $+\infty$,
and so, every point of the Martin boundary is a Busemann point.
\end{proof}
\begin{remark}
Karlsson, Metz, and Noskov~\cite{karlsson}
have shown previously that every boundary point
of a polyhedral normed space is the limit of a geodesic,
and hence a Busemann point. They did this by
characterising the sequences which converge to a boundary point.
\end{remark}
\begin{example}
Consider now $L(x):=\|x\|_\infty^p/p$ with
$\|x\|_\infty:=\max(|x_1|,\cdots,|x_n|)$
and $p>1$.
By Proposition~\ref{prop-carac-bus},
the Busemann points of $(\R^n,\|\cdot\|_\infty)$ with respect
to the basepoint $0$ are of the form:
\[
w:\; x\mapsto 
\min_{i\in I} \epsilon_i(x_i -X_i)+\max_{i\in I} \epsilon_i X_i\enspace,
\]
where $I$ is a non-empty subset of $\{1,\ldots,n\}$, $\epsilon_i=\pm 1$, and the $X_i$ are arbitrary reals. Theorem~\ref{cor-final} shows that any
eigenvector with eigenvalue $\lambda>0$ of the Lax-Oleinik semigroup
can be written as a supremum of maps $c+\vartheta_\lambda w$, where $c\in \R\cup\{-\infty\}$ and $w$ is of the above form. 
For instance, when $n=2$, the functions $w$ are of one of the following forms:
\[
\epsilon_1 x_1,\quad\epsilon_2 x_2,\text{ or } \min(\epsilon_1 (x_1-X_1),\epsilon_2 (x_2-X_2))+\max(\epsilon_1 X_1,\epsilon_2 X_2) \enspace,
\]
with $X_1,X_2\in \R$ and $\epsilon_1=\pm 1,\epsilon_2=\pm 1$.
\end{example}
\begin{remark}
It is natural to ask whether the eigenvectors of the Lax-Oleinik semigroup
$(T^t)_{t\geq 0}$ coincide with the viscosity solutions
of the ergodic Hamilton-Jacobi equation
\[
L^{\star}(\nabla u) =\lambda \enspace , 
\]
where $L^{\star}$ denotes the Legendre-Fenchel transform of $L$.
This is proved in~\cite[Chapter~7]{fathi03} in the different
setting where the space is a compact manifold
and the Lagrangian $L$ can depend on both the position
and the speed but must satisfy certain regularity
and coercivity conditions.
\end{remark}
\bibliographystyle{myalpha}
\bibliography{denumerable}
\end{document}